\documentclass[a4paper,12pt]{book}
\usepackage{csquotes}
\usepackage{url}
 
    \usepackage{amssymb}
  \usepackage{latexsym}
  \usepackage{amsfonts}
  \usepackage{amsthm}
  \usepackage{amsmath}
  \usepackage{mathtools}
\usepackage[]{algorithm2e}
\usepackage{bm}
\usepackage{adjustbox}
\usepackage{changepage}
\usepackage{tablefootnote}
\usepackage{graphicx}
\usepackage{subcaption}
\usepackage{hyperref}
\usepackage{listings}
\newcommand{\siv}{S(\cA; \cP,z)}
\newcommand{\pri}{\mathbb{P}}
\newcommand{\Li}{\text{Li}}
\newcommand{\seta}{\cA_{\leq_x}}
\newcommand{\1}{\textbf{1}}
\newcommand{\al}{\alpha}

\newtheorem{theorem}{Theorem}[section]

\newtheorem{proposition}[theorem]{Proposition}

\newtheorem*{theorem*}{Theorem}
\newtheorem*{corollary*}{Corollary}
\newtheorem*{lemma*}{Lemma}
\newtheorem*{proposition*}{Proposition}

\usepackage{mathtools}

  \newcommand{\cA}{{\mathcal A}}

  \newcommand{\cL}{{\mathcal L}}
  
  \newcommand{\cP}{{\mathcal P}}

  \newcommand{\cS}{{\mathcal S}}

  \newcommand{\R}{{\mathbb R}}

   \newcommand{\N}{{\mathbb N}}

\renewcommand\mod{\text{ mod }}
\newcommand{\lcm}{\text{lcm}}

   \newcommand{\eps}{\epsilon}

  \newcommand{\blanknonumber}{\newpage\thispagestyle{empty}}

  \usepackage{graphics}

\begin{document}

  \frontmatter

  \begin{titlepage}
\begin{center}

\vspace*{\fill} \Huge
                        On a numerical upper bound for the extended Goldbach conjecture
\\
\vfill\vfill\Large
                          David Quarel
\\
\href{mailto:david.quarel@anu.edu.au?subject=Honours Thesis}{\nolinkurl{david.quarel@anu.edu.au}}
\\
\vfill\vfill
                         June 2017 
                          \\
\vfill\vfill \normalsize
         A thesis submitted for the degree of Bachelor of Science (Advanced) (Honours)\\
         of the Australian National University\\
         Revised December 2017
\vfill
         \includegraphics{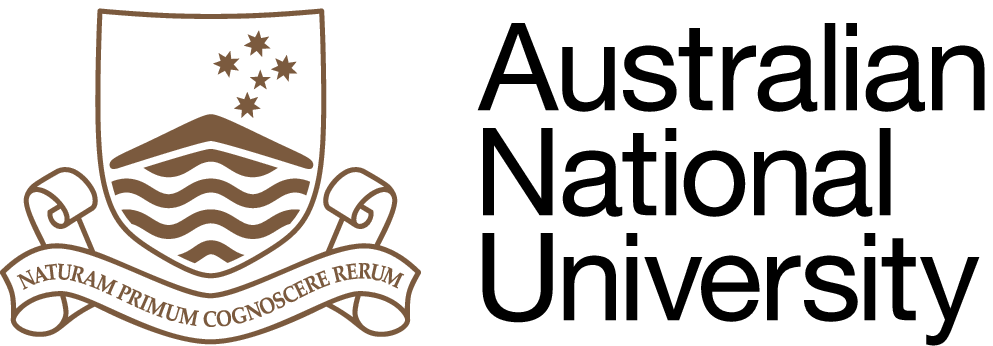}

\end{center}

\end{titlepage}
\blanknonumber

\chapter*{Declaration}\label{declaration}
\thispagestyle{empty}
The work in this thesis is my own except where otherwise stated.

\vspace{1in}

\hfill\hfill\hfill
David Quarel
\hspace*{\fill}
\blanknonumber
  
\chapter*{Acknowledgements}\label{acknowledgements}
\addcontentsline{toc}{chapter}{Acknowledgements}

I wish to thank my supervisor Tim Trudgian. His advice and weekly meetings provided much needed support whenever
I'd get stuck on a problem or become disheartened. I find myself in the most fortunate position of finishing honours
in mathematics with the same person as when I started it, way back in 2013 when I took MATH1115.\footnote{I still remember the hilarious anecdote in class about cutting a cake ``coaxially''.} 

Having the most engaging lecturer during my first year of undergraduate helped to foster my interest in mathematics, and I've yet to see any other lecturer at this university 
take the time to write letters of congratulations to those that did well in their course.

\blanknonumber
  \chapter*{Abstract}\label{abstract}

\addcontentsline{toc}{chapter}{Abstract}

The Goldbach conjecture states that every even number can be decomposed as the sum of two primes.
Let $D(N)$ denote the number of such prime decompositions for an even $N$. It is known
that $D(N)$ can be bounded above by
$$
D(N) \leq C^* \Theta(N), \quad \Theta(N) := \frac{N}{\log^2 N}\prod_{\substack{p|N \\ p>2}} \left( 1 + \frac{1}{p-2} \right)\prod_{p>2} \left( 1 - \frac{1}{(p-1)^2}\right) 
$$
where $C^*$ denotes Chen's constant. It is conjectured \cite{hardy1923some} that $C^*=2$. In 2004, Wu \cite{wu2004chen} showed that $C^* \leq 7.8209$.
We attempted to replicate his work in computing Chen's constant, and
in doing so we provide an improved approximation of the Buchstab function $\omega(u)$,
\begin{align*}
\omega(u) = 1/u,& \quad (1 \leq u \leq 2), \\ 
(u \omega(u))' = \omega(u-1), & \quad  (u \geq 2).
\end{align*}
based on work done by Cheer and Goldston \cite{cheer1990differential}. For each interval $[j,j+1]$,
they expressed $\omega(u)$ as a Taylor expansion about $u=j+1$. We expanded about the point $u=j+0.5$,
so $\omega(u)$ was never evaluated more than $0.5$ away from the center of the Taylor expansion,
which gave much stronger error bounds.  

Issues arose while using this Taylor expansion to compute the required integrals for Chen's constant,
so we proceeded with solving the above differential equation to obtain $\omega(u)$, and then integrating the result. 
Although the values that were obtained undershot Wu's results, 
we pressed on and refined Wu's work by discretising his  
integrals with finer granularity. 
The improvements to Chen's constant were negligible (as predicted by Wu).
This provides experimental evidence, but not a proof, that were Wu's integrals computed on 
smaller intervals in exact form, 
the improvement to Chen's constant would be similarly negligible. 
Thus, any substantial improvement on Chen's constant likely requires
a radically different method to what Wu provided.

\blanknonumber
  \tableofcontents\blanknonumber

\chapter{Notation and terminology}\label{notation}

\renewcommand{\thefootnote}{\fnsymbol{footnote}}

In the following, $p,q$ usually denote primes.\\

\noindent\textbf{Notation}


\newcommand{\nttn}[2]{\item[{\ \makebox[3.18cm][l]{#1}}]{#2}}
\begin{list}{}{ \setlength{\leftmargin}{3.4cm}
                \setlength{\labelwidth}{3.4cm}}

\nttn{GRH}{Generalised Riemann hypothesis.}
\nttn{TGC}{Ternary Goldbach conjecture.}
\nttn{GC}{Goldbach conjecture.}
\nttn{FTC}{Fundamental theorem of calculus.}
\nttn{$f(x)\ll g(x)$}{$g$ is an asymptotic upper bound for $f$, that is, there exists constants $M$ and $c$ such that for all  $x > c$, $|f(x)| \leq M|g(x)|$.}
\nttn{$f(x)\gg g(x)$}{$g$ is an asymptotic lower bound for $f$, that is, there exists constants $M$ and $c$ such that for all $ x > c$, $|f(x)| \geq M|g(x)|$. Equivalent to $g \ll f$.}
\nttn{$f(x)\sim g(x)$}{$g$ is an asymptotic tight bound for $f$, that is, $\lim\limits_{n \to \infty} \frac{f(n)}{g(n)} = 1$. Implies that $f \ll g$ and $g \ll f$.}
\nttn{$C_0$}{Twin primes constant, \S \ref{eq:twinprimes}}
\nttn{$\pi(x)$}{Number of primes $p$ less than $x$.}
\nttn{$\mathbb{P}$}{Set of all primes $p$.}
\nttn{$\mathbb{P}_{\leq n}$}{Set of all numbers with at most $n$ prime factors. Also denoted $\pri_n$ for brevity.}
\nttn{$\siv$}{Sifting function, \S \ref{eq:sieve}.}
\nttn{$\cA$}{Set of numbers to be sifted by a sieve $\siv$.}
\nttn{$\cA_{\leq x}$}{$\{a : a \in \cA, a \leq x \}.$}
\nttn{$\cA_d$}{$\{ a : a \in \seta, a \equiv 0 \mod d \}$ for some square-free $d$.}
\nttn{$\1_E(x)$}{Characteristic function on a set $E$. Equals 1 if $x \in E$, 0 otherwise.}
\nttn{$\cP$}{Usually denotes a set of primes.}
\nttn{$X$}{Shorthand for $X(x)$, denotes the main term in an approximation to a sieve.}
\nttn{$R$}{The remainder, or error term, for an approximation to a sieve.}
\nttn{$C^*$}{Chen's constant, \S\ref{eq:chens_constant}.}
\nttn{$\omega(u)$}{Buchstab's function, \S\ref{eq:buchstab}.}
\nttn{$\gamma$}{Euler--Mascheroni constant.}
\nttn{$\mu(x)$}{M\"{o}bius function.}
\nttn{$\varphi(n)$}{Euler totient function.}
\nttn{$K$}{Linnik--Goldbach constant, \S\ref{eq:linnik}.}
\nttn{$\#A$}{The cardinality of the set $A$.}
\nttn{$(a,b)$}{Greatest common divisor of $a$ and $b$. Sometimes denoted $\gcd(a,b)$ if $(a,b)$ is ambiguous and could refer to an ordered pair.}
\nttn{$\lcm(a,b)$}{Least common multiple of $a$ and $b$.}
\nttn{$f*g$}{Dirichlet convolution of $f$ and $g$, defined as $(f*g)(x) = \sum_{d|n}f(d)g(n/d)$}
\nttn{$a \mid b$}{$a$ divides $b$, i.e.\ there exists some $n$ such that $b = na$.}
\nttn{$a\nmid b$}{$a$ does not divide $b$}
\nttn{$\Li(x)$}{Logarithmic integral, $\Li(x) = \int_2^x \frac{1}{\log t} \; dt$}
\end{list}

%
%
%
%
\blanknonumber

  \mainmatter

  \newcommand{\K}{C^*}

\chapter{History of Goldbach problems} 
  
\section{Origin}
One of the oldest and most difficult unsolved problems in mathematics is the
Goldbach conjecture (sometimes called the Strong Goldbach conjecture), which originated during
correspondence between Christian Goldbach
and Leonard Euler in 1742. The Goldbach conjecture states that ``Every even integer $N>2$ can be expressed as the sum of two primes" \cite{guy1994unsolved}.\\
For example,
\begin{equation*}
4 = 2 + 2 \qquad 8 = 5 + 3 \qquad 10 = 7 + 3 = 5 + 5.
\end{equation*}
Evidently, this decomposition need not be unique.
The Goldbach Conjecture has been verified by computer search for all even $N \leq 4\times 10^{18}$ \cite{goldbachverify}, but 
the proof for all even integers remains open.\footnote{Weakening Goldbach to be true for only sufficiently large $N$ is also an open problem.}
Some of the biggest steps towards solving the Goldbach conjecture
include Chen's theorem \cite{ross1975chen} and the Ternary Goldbach conjecture \cite{helfgott2013ternary}.

\section{Properties of $D(N)$}
 
As the Goldbach conjecture has been unsolved for over 250 years, a lot of work has gone
into solving weaker statements, usually by weakening the statement that $N$ is the sum of
two primes. Stronger statements have also been explored. We define $D(N)$ as the number of ways $N$
can be decomposed into the sum of two primes. The Goldbach conjecture is then equivalent to  $D(2N) \geq 1$ for all $N$.
As $N$ grows larger, more primes exist beneath $N$ (approximately $\frac{N}{\log N}$ many, from the prime number theorem) that could be used to construct a sum for $N$.
Thus, we expect large numbers to have many prime decompositions. And indeed empirically this seems to be the case.

\begin{figure}[ht!]
  \centering
    \includegraphics[width=1\linewidth]{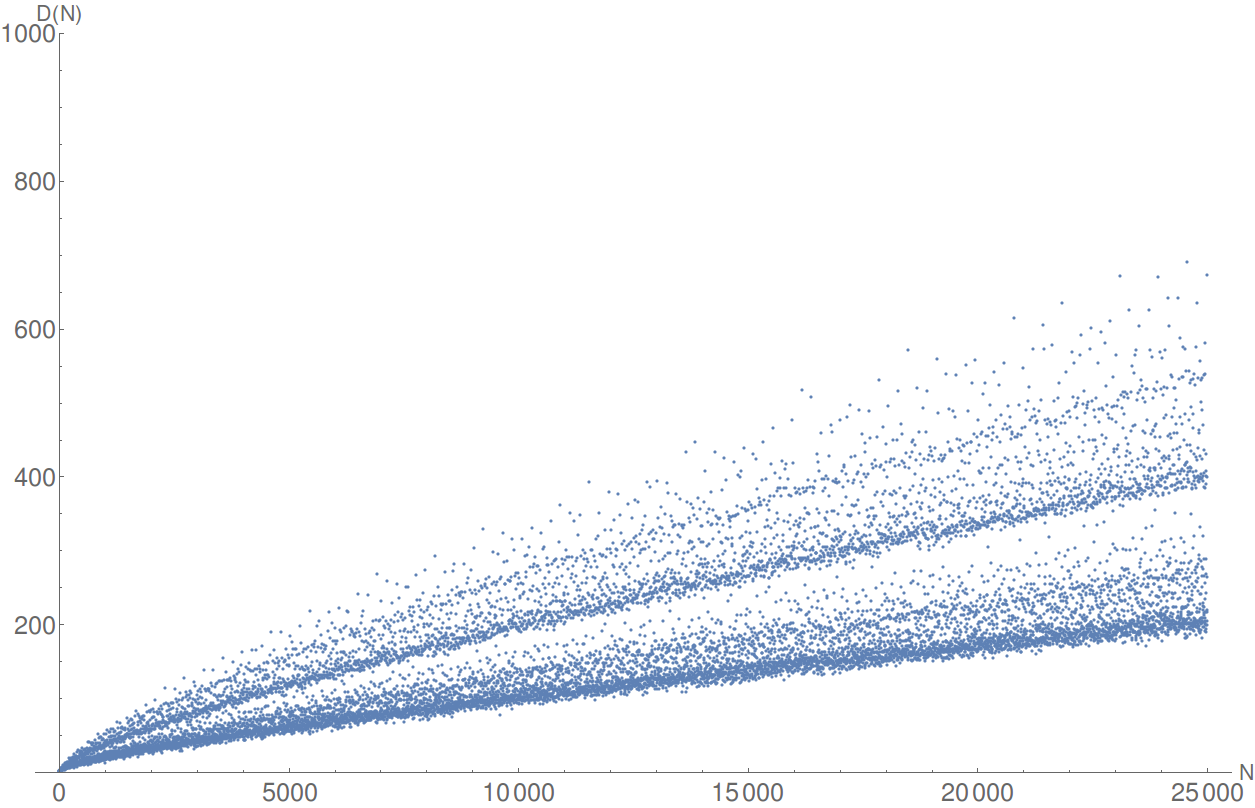}
    \caption{Plot of $D(N)$ for $N \leq 25000$.}
    \label{fig:d(n)v2}
      \vspace*{\floatsep}
  \includegraphics[width=1\linewidth]{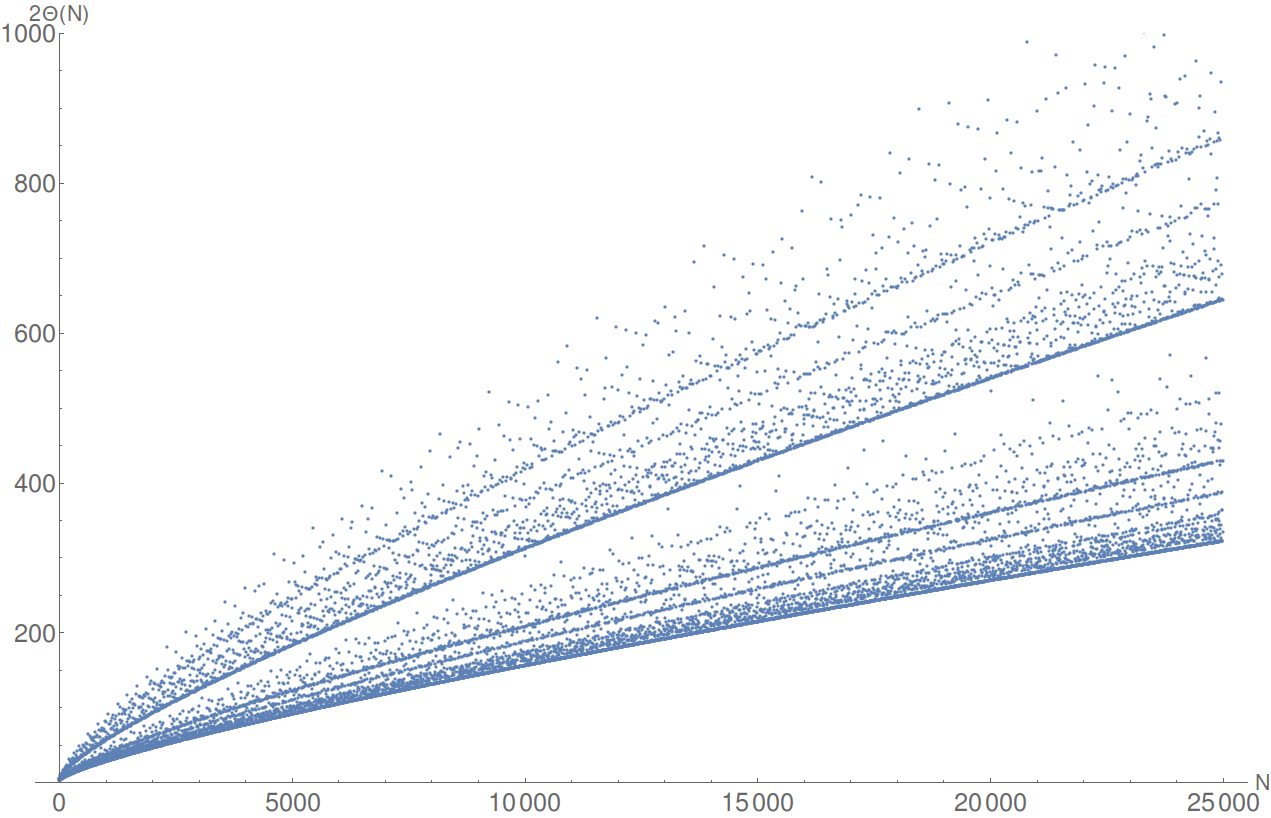}
  \caption{Plot of $2\Theta(N)$, the conjectured tight bound to $D(N)$, for $N \leq 25000$.}
  \label{fig:2Theta(n)_new}
\end{figure}

Below are some examples of prime decompositions of even numbers.
\begin{align*}
36 &= 31 + 5  	& 66 &= 61 + 5 & 90 &= 83 + 7  &= 61 + 29 \\
   &= 29 + 7  	& 	 &= 59 + 7 &    &= 79 + 11 &= 59 + 31 \\
   &= 23 + 13  &     &= 53 + 13 &   &= 73 + 17  &= 53 + 37\\
   &= 19 + 17  &     &=  47 + 19  &  &= 71 + 13 &= 47 + 43 \\
   &          &      &= 43 + 23  &   &= 67 + 23 &         \\
   &          &      &= 37 + 29  &  & &
\end{align*}

Intuitively, we would expect that as $N$ grows large, $D(N)$ should likewise grow large.
The extended Goldbach conjecture \cite{hardy1923some} asks how $D(N)$ grows asymptotically.
As expected, the conjectured formula for $D(N)$ grows without bound as $N$ increases (Figure \ref{fig:d(n)v2}).

\begin{figure}[ht!]
   \centering
     \includegraphics[width=1\linewidth]{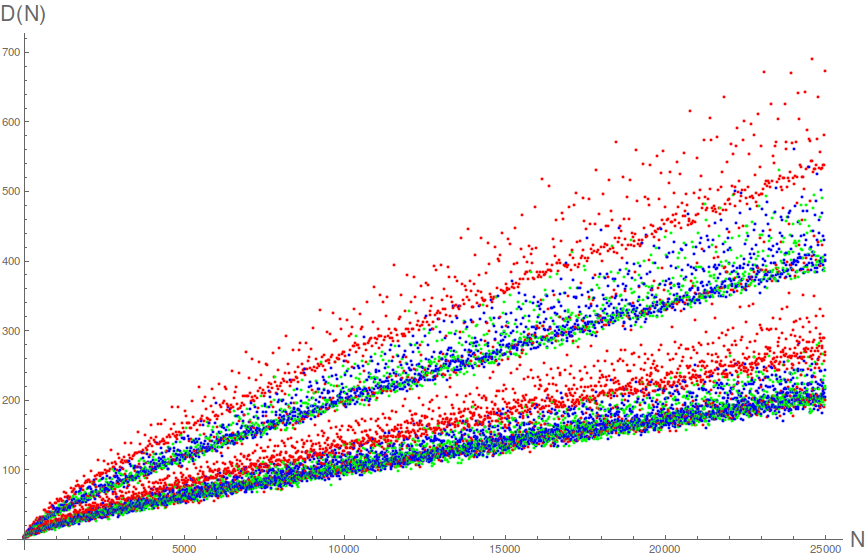}
     \caption{Plot of $D(N)$ for $N \leq 25000$, points coloured based on $N/2 \mod 3$.}
     \label{fig:d(n)_colour}
      \vspace*{\floatsep}
 	 \includegraphics[width=1\linewidth]{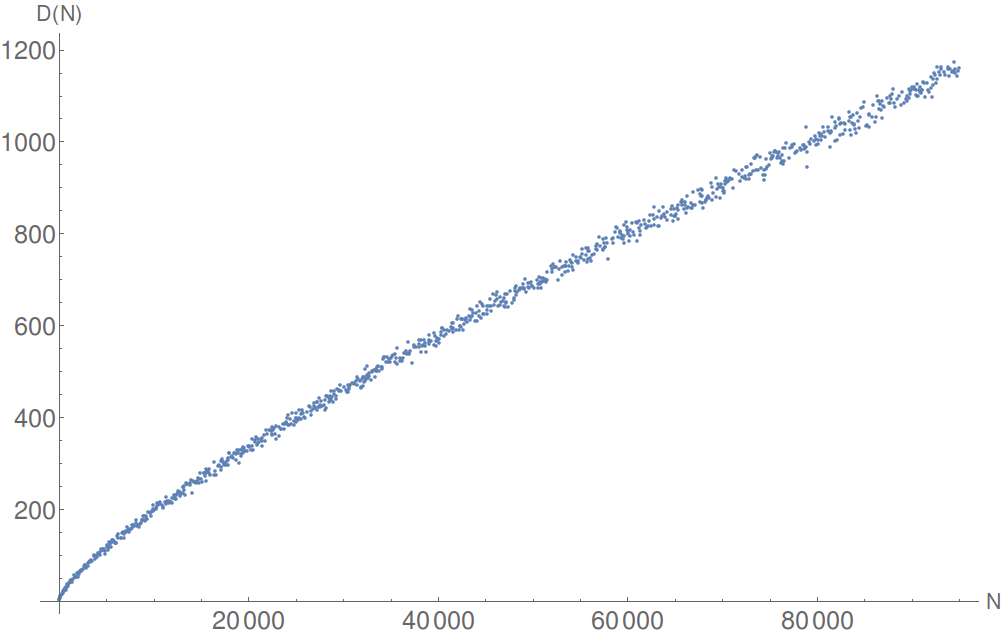}
 	     \caption{Plot of $D(N)$ restricted to $N = 12p, 12p \leq 10^5$}
 	     \label{fig:prime12}
\end{figure}

We can see from the plots that the points of $D(N)$ cluster into bands, and $2\Theta(N)$ also shares this property (Figure \ref{fig:2Theta(n)_new}).
The plot of $D(N)$ called ``Goldbach's Comet'' \cite{fliegel1989goldbach}, which has many interesting properties.
We plot $D(N)$, and colour each point either red, green or blue if $(N/2 \mod 3)$ is 0, 1 or 2 respectively (Figure \ref{fig:d(n)_colour}). \cite{fliegel1989goldbach}.
We can see the bands of colours are visibly separated.
If $D(N)$ is plotted only for prime multiples of an even number \cite{baker2007excel}, say 
$$
N = 12\times 2, 12 \times 3, 12 \times 5, \ldots
$$
the resulting plot is almost a straight line (Figure \ref{fig:prime12}).

\section{Ternary Goldbach conjecture}
The Ternary Goldbach conjecture (TGC) states that ``every odd integer greater than
5 can be expressed as the sum of three primes". This statement is directly implied by the Goldbach conjecture, as 
if every even $n>2$ can be expressed as $n=p+q$, where $p,q$ prime,
we can express all odd integers $n'$ greater than 5 as three primes
$$
n' = n + 3 = p + q + 3.
$$
In 1923, Hardy and Littlewood \cite{hardy1923some} showed that, under the assumption of the generalised Riemann hypothesis (GRH), 
the TGC is true for all sufficiently large odd numbers. They also specified how $D_3(n)$ (the number
of ways $n$ can be decomposed as the sum of three primes) grows asymptotically
\begin{equation}
\label{eq:ternary}
D_3(n) \sim  C_3 \frac{n^2}{\log^3 n}
\prod_{\substack{p|n \\ p\text{ prime}}} \left( 1 - \frac{1}{p^2 - 3p + 3}\right),
\end{equation}
where
$$
C_3 = \prod_{\substack{p > 2 \\ p \text{ prime}}} \left(1 + \frac{1}{(p-1)^3} \right).
$$
Using the Taylor expansion of $\exp(x)$ and the fact that every term in the sum is positive, we obtain the inequality
$$
1+x \leq 1 + x + \frac{x^2}{2!} + \ldots = \sum_{n=0}^{\infty} \frac{x^n}{n!} = \exp(x).
$$
Hence by induction, for all non-negative sequences $\{a_1,\ldots,a_n\}$
\begin{equation}
\label{eq:prod-ineq}
\prod\limits_{i=1}^{n} (1 + a_i) \leq \exp \left( \sum_{i=1}^n a_i\right).
\end{equation}
Therefore we may note that $C_3$ must converge, as
$$
C_3 < \prod_{n \geq 2 } \left(1 + \frac{1}{n^3} \right) \leq \exp \left( \sum_{n=2}^\infty \frac{1}{n^3}\right)
< \infty.
$$
We also note that the product over all primes dividing $n$ in (\ref{eq:ternary})  grows slowly compared to the main
term, as for all $n \geq 2$
\begin{align*}
n^2 - 3n + 3 &\geq \frac{1}{2}n \\
1- \frac{1}{n^2 - 3n + 3} &\geq 1-\frac{2}{n}\\
\prod_{\substack{p|n \\ p \text{ prime}}} \left( 1 - \frac{1}{p^2 - 3p + 3}\right)
\geq
\prod_{\substack{p|n \\ p \text{ prime}}} \left( 1 - \frac{2}{p}\right)
&\geq \prod_{\substack{2 < p \leq n \\ p \text{ prime}}} \left( 1 - \frac{2}{p}\right).
\end{align*}
We have from Mertens' theorem \cite{de2012analytic} that
\begin{equation}
\label{eq:mertens}
\prod_{\substack{p \leq n \\ p \text{ prime}}} \left( 1 - \frac{1}{p} \right)
\sim \frac{e^{-\gamma}}{\log n}.
\end{equation}
Combining this with the following inequality
\begin{equation}
\label{eq:poly_inequality}
1-\frac{2}{x} \geq \frac{1}{2}(1-\frac{1}{x})^2
\text{ for } x \geq 3,
\end{equation}
we obtain
$$
\prod_{\substack{2 < p \leq n \\ p \text{ prime}}} \left( 1 - \frac{2}{p}\right)
\geq 
\frac{1}{2}\prod_{\substack{2 < p \leq n \\ p \text{ prime}}} \left( 1 - \frac{1}{p}\right)^2
\sim \frac{e^{-2 \gamma}}{2 \log ^2 n },
$$
which gives an asymptotic lower bound on $D_3(n)$
$$
D_3(n) \gg  \frac{n^2}{\log^5 n}.
$$
This gives a stronger version of the TGC, as it implies that the number of different ways that an odd integer 
can be decomposed as the sum
of three primes can be arbitrarily large. 
In 1937, Vinogradov improved this result by removing the dependency on GRH \cite{vinogradov1937}.

An explicit value for how large $N$ needs to be before the TGC holds was found by Borozdin, who showed that $N > 3^{3^{15}}$ is sufficient \cite{golomb1985invincible}. 
In 1997, the TGC was shown to be conditionally true for all odd $N > 5$ by 
Deshouillers et al, \cite{deshouillers1997complete} under the assumption of GRH.
The unconditional TGC was further improved by Liu and Wang \cite{M2002} in 2002, who showed that Borozdin's constant 
can be reduced to $N > e^{3100}$. Thus, in principle, all that was needed to prove the TGC was to
check that it held for all $N < e^{3100} \approx 10^{7000}$. Given that the number of particles in the universe 
$\approx 10^{80}$, no improvements to computational power would likely help. Further mathematical work 
was necessary.

In 2013, Helfgott reduced the constant sufficiently as to allow computer
assisted search to fill the gap, obtaining $N > 10^{27}$ \cite{helfgott2013ternary}.
When combined with his work numerically verifying the TGC for all 
$N<8.875 \times 10^{30}$, \cite{helfgott2013ternary}
Helfgott succeeded in proving the TGC.

\section{Extended Goldbach conjecture}
Hardy and Littlewood's paper \cite{hardy1923some} discuss how one can directly obtain an asymptotic formula 
(conditional on GRH) for the numbers of ways even $N$ can be decomposed into the
sum of four primes. Let $D_4(N)$ be the number of such ways.
\begin{equation}
D_4(N) \sim \frac{1}{3}C_4 \frac{n^3}{\log^4 n} \prod_{\substack{p | n \\ p>2}} 
	\left( 1 + \frac{1}{(p-2)(p^2-2p+2)} \right),
\end{equation}
where 
$$
C_4 = \prod_{p>2} \left(1 - \frac{1}{(p-1)^4} \right).
$$
The asymptotic formulae for $D_3(N)$ and $D_4(N)$ are very similar 
(checking convergence and the long term behaviour of $D_4(N)$ is the same proof as for $D_3(N)$). Hardy--Littlewood 
claim that one can easily generalise their work to obtain an asymptotic formula for $D_r(N)$ for $r>2$. However,
for $r=2$, the results do not easily generalise:

\begin{displayquote}
``It does not fail \emph{in principle}, for
it leads to a definite result which appears to be correct; but we cannot overcome 
the difficulties of the proof..."\cite[p. 32]{hardy1923some}
\end{displayquote}
The following asymptotic formula for $D(N)$ was thus conjectured \cite{hardy1923some}
\begin{equation}
\label{eq:hardyconj}
D(N) \sim 2 \Theta(N),
\end{equation}
where
$$
\Theta(N) := \frac{C_N N}{\log^2 N},
\quad \quad
C_N := C_0 \prod_{\substack{p|N \\ p>2}} \left( 1 + \frac{1}{p-2} \right),
\quad C_0 := \prod_{p>2} \left( 1 - \frac{1}{(p-1)^2}\right).
$$
Now $C_0$ denotes the twin primes constant (see \ref{eq:twinprimes}), and 
each term in the infinite product in $C_N$ is greater than 1. Thus, $C_N > C_0$,
which provides an asymptotic lower bound of
\begin{equation}
\label{eq:lower_bound_D(N)}
D(N) \gg \frac{N}{\log^2 N},
\end{equation}
Hence, (\ref{eq:hardyconj}) implies the Goldbach conjecture (for large $N$),
which should indicate the difficulty of the problem.

There has been progress in using (\ref{eq:hardyconj}) as a way to construct an upper bound for $D(N)$,
by seeking the smallest values of $C^*$ (hereafter referred to as Chen's constant) such that
\begin{equation}
\label{eq:chens_constant}
D(N) \leq \K \Theta(N)
\footnote[3]{Some authors instead include the factor of 2 inside $C^*$, so they instead 
	assert $D(N) \sim \Theta(N)$ for (\ref{eq:hardyconj}), and
	$D(N) \leq 2 C^* \Theta(N)$ for (\ref{eq:chens_constant}). }
\end{equation}
Upper bounds for $\K$ have been improved over the years, but recent improvements 
have been small\footnote[2]{We remind the reader than $\K$ is conjectured to be 2.}
(see Table \ref{tab:K}). 
Recent values have been obtained by various sieve theory methods 
by constructing large sets of weighted inequalities (the one in Wu's paper
has 21 terms!). However, the reductions in $\K$ are minuscule. One would suspect that any further improvements
would require a radically new method, instead of complicated inequalities with more terms.

\begin{table}[ht!]
	\centering
    \begin{tabular}{ | l | l  | l |}
    \hline
    $\K$ & Year & Author \\ \hline
    $16 + \epsilon$ & 1949 & Selberg \cite{selberg1952} \\
    12 & 1964 & Pan \cite{pan1964new} \\
    8 & 1966 & Bombieri \& Davenport \cite{bombieri1966small} \\
    7.8342 & 1978 & Chen \cite{chen1978goldbach} \\
    7.8209 & 2004 & Wu \cite{wu2004chen} \\
    \hline
    \end{tabular}
    \caption{Improvements of $\K$ in (\ref{eq:chens_constant}) over time, conjectured that $C^* = 2$. \label{tab:K}}
\end{table} 

Wu Dong Hua \cite{dongwrong} claimed to have $C^* \leq 7.81565$, but some functions he defined to compute $C^*$
failed to have some necessary properties \cite{tim_email}. 

\section{The Linnik--Goldbach problem}

The Linnik--Goldbach problem \cite{platt2015linnik} is another weaker form of the GC, which asks 
for the smallest values of $K$
such that for all sufficient large even $N$,
\begin{equation}
\label{eq:linnik}
N = p_1 + p_2 + 2^{e_1} + 2^{e_2} + \ldots + 2^{e_r}, \; r \leq K.
\end{equation}
That is, we can express $N$ as the sum of two primes and at most $K$ powers of two
(we refer to $K$ as the Linnik--Goldbach constant).
In 1953, Linnik \cite{linnik1953addition} proved that there exists some finite $K$ such that the statement holds, 
but did not include an explicit value for $K$. (See Table \ref{tab:linnik} for historical improvements on $K$.)
The methods of Heath-Brown and Puchta \cite{heath2002integers}, and of Pintz and Ruzsa \cite{pintz2003linnik} show that $K$ satisfies (\ref{eq:linnik})
if
$$
\lambda^{K-2} < \frac{C_3}{( \K - 2)C_2 + cC_0^{-1}\log 2},
$$
for particular 
constants 
 $\lambda, C_0,\K,C_2,C_3,c$.
 
 Here, $C_0$ is the twin primes constant, given by
 \begin{equation}
 \label{eq:twinprimes}
 C_0 = \prod_{\substack{p>2 \\ p \text{ prime}}} \left( 1 - \frac{1}{(p-1)^2}\right)
 \approx 0.66016 \ldots
 \end{equation}
 The infinite product for $C_0$ is easily verified as convergent, as
 for all integers $n > 2$
 $$
 0 < 1-\frac{1}{(n-1)^2} < 1.
 $$ 
 So $C_0$ is an infinite product, with each term strictly between 0 and 1. Therefore,
 $0 \leq C_0 < 1$.
 The value of $C_0$ is easily computed: Wrench provides $C_0$ 
 truncated to 42 decimal places \cite{wrench1961evaluation}.
 As far as lowering the value of $K$, improvements on the other constants
 will be required.
 
 Now, $C_2$ is defined by 
 $$
 C_2 = \sum_{d=0}^{\infty} \left( \frac{|\mu(2d+1)|}{\epsilon(2d+1)} \prod_{\substack{p|d \\ p>2}} \frac{1}{p-2} \right),
 $$
 where
 $$
 \eps(n) = \min_{v \in \N} \{ v : 2^v \equiv 1\mod d \}.
 $$
 and $\mu(x)$ is the M\"{o}bius function.
 It has been shown by Khalfalah--Pintz \cite{pintz2004goldbach}
 that
 \begin{equation}
 \label{eq:pintz}
 1.2783521041 < C_2C_0 < 1.2784421041,
 \end{equation}
 So since $C_0$ is easily computed, we can convert (\ref{eq:pintz}) to
 \begin{equation}
 1.93642 < C_2 < 1.93656.
 \end{equation}
 
 It is remarked by Khalfalah--Pintz that estimating $C_2$ is very difficult, and that any further progress is unlikely.
 The remaining constants are defined by similarly complicated expressions (see \cite{platt2015linnik}).
 Platt and Trudgian \cite{platt2015linnik} make improvements on both $C_3$ and $\lambda$,
 by showing that one can take
 $$
 (C_3,\lambda) = (3.02858417,0.8594000)
 $$
 and obtain unconditionally that $K \geq 11.0953$, a near miss for $K = 11$. Platt--Trudgian
 remark that any further improvements in estimating $C_3$ using their method with more computational power
 would be limited.
 Obtaining $K = 11$ by improving Chen's constant (assuming all others constants used by Platt--Trudgian are the same) 
 would require $C^* \leq 7.73196$, which is close to the best known value of $C^* \leq 7.8209$. 
 This provides a motivation for reducing $C^*$.
 
\begin{table}[ht!]
	\centering
    \begin{tabular}{ | l | l  | l | l |}
    \hline
    $K$ 		& $K$ assuming GRH & Year & Author \\ \hline
    - & $<\infty$ 				 & 1951 & Linnik \cite{linnik1951prime} \\
    $<\infty$ & - 				 & 1953 & Linnik \cite{linnik1953addition} \\
    54000 & - & 1998 & Liu, Liu and Wang \cite{liu1998number} \\
    25000 & - & 2000 & H.Z.Li \cite{li2000number} \\
        - & 200 & 1999 & Liu, Liu and Wang \cite{liu1999almost} \\
    2250 & 160 & 1999 & T.Z.Wang \cite{wang1999linnik} \\
    1906 & - & 2001 & H.Z.Li \cite{li2001number} \\
    13 & 7 & 2002 & Heath-Brown and Puchta$^*$ \cite{heath2002integers} \\
    12 & - & 2011 & Liu and L\"{u} \cite{liu2011density} \\
    \hline
    \end{tabular}
    \caption{Values of $K$ in (\ref{eq:linnik}) over the years\label{tab:linnik}}
\end{table} 

\newcommand\blfootnote[1]{%
  \begingroup
  \renewcommand\thefootnote{}\footnote{#1}%
  \addtocounter{footnote}{-1}%
  \endgroup
}

\section{History of Romanov's constant}
In 1849, de Polignac conjectured that every odd $n>3$ can be written as a sum of an \emph{odd} prime and a power of two.
He found a counterexample $n=959$, but an interesting question that follows is how many counterexamples are
there? Are there infinitely many? If so, how common are they among the integers? If we let
\begin{equation}
A(N) = \#\{ n : n<N,  n=2^m + p \}
\end{equation}
\begin{equation}
\label{eq:romanov}
R_- = \liminf_{N \to \infty} \frac{A(N)}{N}
\quad
R = \lim_{N \to \infty} \frac{A(N)}{N}
\quad
R_+ = \limsup_{N \to \infty} \frac{A(N)}{N}
\footnote{Some authors denote $R = \lim_{n \to \infty}2A(N)/N$ (and similar for $R_+$ and $R_-$), in which case the trivial bound is $R_+ \leq 1$.}
\end{equation}
then $R$ denotes the density of these decomposable numbers, and $R_+$ and $R_-$ provide upper and lower
bounds respectively.

\blfootnote{$\!\!\!\!\!\!^*$This was also proved independently by Pintz and Ruzsa \cite{pintz2003linnik} in 2003.
            They also claim to have $K \leq 8$, but this is yet to appear in literature.}

Romanov \cite{romanoff1934einige} proved in 1934 that $R_->0$, though
he did not provide an estimate on the value of $R$. 
An explicit lower bound on $R$ was not shown until 2004, with 
Chen--Sun \cite{chen2004romanoff} who proved $R_- \geq 0.0868$. 
Habsieger--Roblot \cite{habsieger2006integers}
improved both bounds, obtaining $0.0933 < R_- \leq R_+ < 0.49045$.
Pintz \cite{pintz2006romanov} obtained $R_- \geq 0.09368$, but under the assumption of Wu Dong Hua's 
value of Chen's constant $C^* \leq 7.81565$, the proof of which is flawed \cite{tim_email}. 
 
Since $2n \neq p + 2^k$ for an odd prime $p$, it is trivial to show that $R_+ \leq 1/2$.
Van de Corput \cite{corput1950} and Erd\H{o}s \cite{erdos1950} proved in 1950 that $R_+ < 0.5$.
In 1983, Romani \cite{romani1983computations} computed $A(N)/N$ for all $N < 10^7$ and noted that
the minima were located when $N$ was a power of two. He then constructed an approximation to $A(N)$,
by using an approximation to the prime counting function $\pi(x)$
$$
\pi(x) \sim \text{Li}(x) := \int_2^x \frac{1}{\log t} \; dt = \sum_{m=1}^{\infty} \frac{(m-1)!}{\log^m x} x,
$$
then using the Taylor expansion of $\text{Li}(x)$ to obtain a formula for $A(N)$ with unknown coefficients $b_i$
$$
A(N) = R x + \sum_{i=1}^\infty b_i \frac{x}{(\log x)^i}.
$$
By using the precomputed values for $A(N)$, the unknown coefficients may be estimated, 
thus extrapolating values for $A(N)$. Using this method, Romani conjectured that 
$R \approx0.434$.

Bomberi's \cite{romani1983computations} probabilistic approach, 
which randomly generates probable primes
(as they are computationally faster to find than primes,
and pseudoprimes\footnote{A composite number that a probabilistic primality test incorrectly asserts is prime.}
are uncommon enough that it does not severely alter the results)  
in the interval $[1,n]$,
in such a way that the distribution closely matches that of the primes, is also used to compute $R$. 
The values obtained closely match Romani's work, and so it is concluded that $R \approx 0.434$ is a reasonable estimate.

The lower bound to Romonov's constant $R_-$, the Linnik--Goldbach constant $K$ and Chen's constant $C^*$ are all related.
There is a very complicated connection relating $R_-$ and $K$ \cite{pintz2006romanov}, but here we prove a simple
property relating $R_-$ and $K$.
\begin{proposition*}
If $R_- > 1/4$, then $K \leq 2$ \cite{pintz2006romanov}. 
\begin{proof}
If for an even $N$ we have that more than $1/4$ of the numbers up to $N$
can be written as the sum of an odd prime and a power of two, then
no even $n$ can be written as $n = p + 2^k$ since $p$ is odd. So more than $1/2$ of
the odd numbers up to $N$ can be written in this way. If we construct a set of pairs of
odd integers 
$$
\left\{(1,N-1), (3,N-3), \ldots, \left(\frac{N-r}{2},\frac{N+r}{2} \right) \right\},
$$
where $r = N \mod 4$, then by the pigeon-hole principle, there must exist a pair $(k,N-k)$
such that
$$
k = p_1 + 2^{v_1}, \quad N-k = p_2 + 2^{v_2}.
$$
Therefore
$$
N = k + (N-k) = p_1 + p_2 + 2^{v_1} + 2^{v_2}. \qedhere
$$
\end{proof}
\end{proposition*}

Pintz also proved the following theorems \cite{pintz2006romanov}.
\begin{theorem}
If $C^* = 2$ holds, then $R_- \geq 0.310987$.
\end{theorem}
This can be combined with the proposition above to obtain a relationship between $C^*$ and $K$,
although Pintz also shows a weaker estimate on $C^*$ is sufficient.
\begin{theorem}
If $C^* \leq 2.603226$, then $R_- > 1/4$ and consequently $K \leq 2$.
\end{theorem}
As an explicit lower bound for $R$ did not appear until 2004, and this bound is much lower
than the expected $0.434$ given by Romani, attempting to improve $K$
by improving $R_-$ appears a difficult task. 

\section{Chen's theorem}
Chen's theorem \cite{chen1973}
is another weaker form of Goldbach, which states that every sufficiently large even integer is the sum of two primes, or a prime and a semiprime (product of two primes).
By letting $R(N)$ denotes the number of ways $N$ can be decomposed in this manner, Chen attempts to get a lower bound on the 
number of ways $N$ can be decomposed as $N = p + P_3$, where $P_3$ is some number with no more than three prime factors, using sieve theory methods.
He then removes the representations where $N = p + p_1p_2p_3$ \cite{halberstam1974sieve}. 
By doing so, Chen obtains that there exists some constant $N_0$ such that if $N$ is even and
$N>N_0$, then
$$
\#\{p : p \leq N, n-p \in \pri_2\} > 0.67 \Theta(N),
$$
where $\Theta(N)$ is the same function used for the conjectured tight bound of $D(N)$ (see \ref{eq:hardyconj}).
We have the asymptotic lower bound (see \ref{eq:lower_bound_D(N)})
$$
\Theta(N) \gg \frac{N}{\log^2 N}.
$$
so this implies that for all sufficiently large even $N$
$$
N = p + q \text{ for } q \in \pri_2
$$
In 2015, Yamada \cite{yamada2015explicit} proved that letting $N_0 = \exp \exp 36$ is sufficient for the above to hold.

  \chapter{Preliminaries for Wu's paper}\label{cha:sieve}

\section{Summary}

Wu \cite{wu2004chen} proved that $C^* \leq 7.8209$.
He obtained this result by expressing the Goldbach conjecture as a linear sieve,
and finding an upper bound on the size of the sieved set. 
The upper bound constructed is very complicated, using a series of weighted
inequalities with 21(!) terms. The approximation to the sieve is computed using numerical integration,
weighting the integrals over judicially chosen intervals. Few intervals (9), over which to
discretise the integration are used, thereby reducing the problem to a linear
optimisation with 9 equations.
\begin{quotation}
\noindent
``If we divide the interval $[2,3]$ into more than 9 subintervals, we can certainly obtain
a better result. But the improvement is minuscule.''\cite[p. 253]{wu2004chen}
\end{quotation}
The main objective of this paper is to quantify how ``minuscule'' the improvement is.

\section{Preliminaries of sieve theory}
A sieve is an algorithm that takes a set of integers and \emph{sifts out}, or removes,
particular integers in the set based on some properties.

The classic example is the \emph{Sieve of Eratosthenes} \cite[Chpt 1]{heath2002lectures}, an algorithm devised
by Erathosthenes in the 3rd century BC to generate all prime numbers below
some bound. The algorithm proceeds as follows: Given the first $n$ primes 
$\cP = \{ p_1 ,\ldots, p_n\}$, list the integers from 2 to $(p_n + 2)^2 - 1$.
For each prime $p$ in the set $\cP$, strike out all multiples of $p$ in the list
starting from $2p$. All remaining numbers are primes below $(p_n + 2)^2$, as the
first number $k$ that is composite and not struck out by this procedure must share no
prime factors with $\cP$. So $k = p_{n+1}^2 \geq (p_n + 2)^2$.

This method of generating a large set of prime numbers is efficient, as by striking off
the multiples of each prime from a list of numbers, the only operations used is addition and
memory lookup. Contrast this with primality testing via trial division, as integer division
is slower than integer addition. O'Neill \cite{o2009genuine} showed that to generate all primes
up to a bound $n$, the Sieve of Eratosthenes takes $O(n \log \log n)$ operations 
whereas repeated trial division would take $O(n^\frac{3}{2} /\log^2 n )$ operations.
Moreover, the Sieve of Eratosthenes does not use division in computing the primes,
only addition. 

Sieve theory is concerned with estimating the size of the 
sifted set. Formally, given a finite set $\cA \subset \N$, a set of primes $\cP$ and some number $z \geq 2$,
we define the sieve function \cite{heath2002lectures}
\begin{equation}
\label{eq:sieve}
S(\cA; \cP,z) = \# \{ a \in \cA : \forall p \in \cP \text{ with } p< z, (a,p) = 1 \},
\end{equation}
i.e.\ by taking the set $\cA$ and removing the elements that are divisible by some
prime $p \in \cP$ with $p < z$, the number of elements left over is $\siv$.
Removing the elements of $\cA$ in this manner is called \emph{sifting}. Analogous to
a sieve that sorts through a collection of objects and only lets some through, the 
unsifted numbers are those in $\cA$ that do not have a prime factor in $\cP$.
Many problems in number theory can be re-expressed as a sieve, and thus attacked. 
We provide an example related to the Goldbach conjecture 
included in \cite{heath2002lectures}.
For a given even $N$, let
$$
\cA = \{ n(N-n) : n \in \N, 2 \leq n \leq N-2\}, \qquad \cP = \pri.
$$
Then 
$$
S(\cA,\pri;\sqrt{N}) = \# \{ n(N-n): n \in \N, p < \sqrt{N} \Rightarrow (n(N-n),p) = 1 \}.
$$
Now since $(n(N-n),p)= 1 \Rightarrow (n,p)=1 \text{ and } (N-n,p)=1$ for all $p < \sqrt{N}$,
both $n$ and $N-n$ must be prime, as if not, either $n$ or $N-n$ has a prime factor $\geq \sqrt{N}$,
which implies $n \geq N$ (impossible) or $N-n \geq N$ (also impossible).
Hence
$$
S(\cA,\pri;\sqrt{N}) = \# \{ p \in \pri: (N-p) \in \pri, p < \sqrt{N} \}, 
$$
Now for all primes $p$, if  $p < \sqrt{N}$, then $p \leq N/2$, hence
\begin{align*}
S(\cA,\pri;\sqrt{N}) &\leq  \# \{ p \in \pri: (N-p) \in \pri , p \leq N/2\}, \\
				&\leq  \# \{ p \in \pri: (N-p) \in \pri , p \leq N-p\} = D(N).
\end{align*}
So this particular choice of sieve is a lower bound for $D(N)$.
If we have a particular (possibly infinite) subset of the integers $\cA$,
we may wish to take all numbers in $\cA$ less that or equal to $x$, denoted\footnotemark
$\seta$
$$
\cA_{\leq x} = \{a \in \cA : a \leq x \},
$$
\footnotetext{In the literature, it is common to 
use $\cA$ instead of $\cA_{\leq x}$, and it is implicitly
understood that $\cA$ only includes the numbers less than some $x$.}
and ask how fast $\# \cA_{\leq x}$ grows with $x$.
For example, if $\cA$ were the set of all even numbers, 
$$
\cA = \{ n \in \N : n \text{ even} \} = \{0,2,4,\ldots \},
$$
the size of the set of $\cA_{\leq x}$ would be computed as:
$$
\# \cA_{\leq x} = \# \{a \in \cA : a \leq x \} = \left\lfloor \frac{x}{2} \right\rfloor + 1.
$$
It is usually more difficult to figure out how sets grow. Likewise,
an exact formula is usually impossible. Therefore, the best case is usually an asymptotic
tight bound.
For example, if $\cA= \pri$, the set of all primes, then the prime number
theorem \cite[p. 9]{apostol1998introduction} gives
$$
\# \{ p \in \pri : p \leq n \} = \pi(x) \sim \frac{x}{\log x}.
$$
The main aim of sieve theory is to decompose $\# \cA_{\leq N}$ as the sum of a main term $X$ and an error term $R$,
such that $X$ dominates $R$ for large $N$. We can then obtain an asymptotic formula for 
$\# \cA_{\leq N}$ \cite{halberstam1974sieve}. 
This partially accounts for the many
proofs in number theory that only work for some extremely large $N$, where
$N$ can be gargantuan (see the Ternary Goldbach conjecture in Chapter 1),
as the main term may only dominate the error term for very large values.
In some cases it can only be shown that the main term eventually dominates the error term,
no explicit value as to how large $N$ needs to be before this occurs is required.

Now, for a given subset $\seta$, we choose a main term\footnotemark $X$ that is hopefully a good approximation
of $\# \seta$.

\footnotetext{Note that since $X$ depends on $x$, we would normally write $X(x)$, but we wished to stick
to convention.}
We denote
$$
\cA_d = \{ a : a \in \seta, a \equiv 0 \mod d \},
$$
for some square-free number $d$.
For each prime $p$, we assign a value for the function $\omega(p)$, with the constraint that $0 \leq \omega(p) < p$, so that
$$
\# \cA_p \approx \frac{\omega(p)}{p} X.
$$
By defining 
$$
\omega(1) = 1, \quad \omega(d) = \prod_{p|d} \omega(p),
$$
for all square-free $d$, we ensure that $\omega(d)$ is a multiplicative function.
We now define the error term (sometimes called the remainder term):
$$
R_d = \#\cA_d - \frac{\omega(d)}{d} X.
$$
So we have that
$$
\#\cA_d = \frac{\omega(d)}{d}X + R_d,
$$
where (hopefully) the main term $\frac{\omega(d)}{d}X$ dominates the error term 
$R_d$.
Given a set of prime numbers $\cP$, we define 
$$
P(z) = \prod_{p \in \cP, p <z} p.
$$

We can rewrite the sieve using two identities of the M\"{o}bius function 
and multiplicative functions.
\begin{theorem}\cite[Thm 2.1]{apostol1998introduction}
If $n \geq 1$ then 
\begin{equation}
\label{eq:sum_mobius}
\sum_{d|n} \mu(d) =
\begin{cases}
1 & n=1, \\
0 & n>1. 
\end{cases}
\end{equation}
\begin{proof}
The case $n=1$ is trivial. For $n>1$, write the prime decomposition
of $n$ as $n = p_1^{e_1} \ldots p_n^{e_n}$. Since $\mu(d)=0$ when $d$ is divisible
by a square, we need only consider divisors $d$ of $n$ of the form $d = p_1^{a_1}\ldots p_n^{a_n}$, where each $a_i$ is either zero or one, (i.e., divisors that are products
of distinct primes). Enumerating all possible products made from $\{p_1 ,\ldots, p_n \}$,
we obtain
\begin{align*}
\sum_{d|n} \mu(d) &= \mu(1) + \Big( \mu(p_1) + \ldots + \mu(p_n) \Big) + \Big( \mu(p_1 p_2) \\  \qquad & + \ldots + \mu(p_{n-1}p_n) \Big) + \ldots + \mu(p_1 \ldots p_n) \\
&= 1 + \binom{n}{1}(-1) + \binom{n}{2}(1) + \ldots + \binom{n}{n}(-1)^n \\
& = \sum_{i=0}^n \binom{n}{i} (-1)^i = (1-1)^n = 0.
\end{align*}
 
\end{proof}
\end{theorem}

\begin{theorem}\cite[Thm 2.18]{apostol1998introduction}
If $f$ is multiplicative
\footnote{A function $f$ is multiplicative if $f(m)f(n) = f(mn)$ whenever $\gcd(n,m)=  1$.} 
 we have
\begin{equation}
\label{eq:mobius_prod}
\sum_{d|n} \mu(d) f(d) = \prod_{p|n}(1-f(p)).
\end{equation}
\begin{proof}
Let
$$
g(n) = \sum_{d|n} \mu(d)f(d).
$$
Then $g(n)$ is multiplicative, as
\begin{align*}
g(pq) &= \sum_{d|pq} \mu(d) f(d)
= \mu(1)f(1) + \mu(p)f(p) + \mu(q)f(q) + \mu(pq)f(pq), \\
&= 1 - f(p) - f(q) + f(p)f(q) = (1 - f(p))(1- f(q)), \\
&= \left( \sum_{d|p} \mu(d) f(d) \right) \left( \sum_{d|q} \mu(d) f(d) \right)
=g(p)g(q).
\end{align*}
So by letting $n = p_1^{a_1} \ldots p_n^{a_n}$, we have
$$
g(n) = \prod_{i=1}^n g(p_i^{a_i})
= \prod_{i=1}^n \sum_{d |p_i^{a_i}} \mu(d) f(d).
$$
Now since $\mu(d)$ is zero for $d = p^2, p^3,\ldots$, the only non-zero terms that are
divisors of $p_i^{a_i}$ are $1$ and $p_i$.
$$
g(n) = \prod_{i=1}^n (\mu(1) f(1) + \mu(p_i)f(p_i))
= \prod_{i=1}^n (1 -f(p_i)) = \prod_{p|n} (1-f(p)).
$$
\end{proof}
\end{theorem}

Now given this, we can rewrite the sieve as
\begin{equation}
\label{eq:rewrite_sieve}
\cS(\cA,\cP,x) = \sum_{\substack{n \in \cA \\ (n,P(z)) = 1}} 1
= \sum_{n \in \cA} \sum_{\substack{ d | n \\ d | P(z)}} \mu(d) 
= \sum_{d|P(z)} \mu(d) \# \cA_d,
\end{equation}
as we sum over all $n \in \cA$ such that $a$ is square-free, as
the M\"{o}bius function is only non-zero for square-free inputs.
Using the approximation to $\# \cA_d$, we obtain
\begin{equation}
\label{eq:sieve_main_error}
\begin{split}
\cS(\cA,\cP,x) &= X \sum_{d|P(z)} \frac{ \mu(d)  \omega(d)}{d} + \sum_{d|P(z)} \mu(d) R_d, \\
 &= X \prod_{p|P(z)} \left(1 - \frac{\omega(p)}{p} \right) + \sum_{d|P(z)} \mu(d) R_d.
\end{split}
\end{equation}
Writing 
$$
W(z;\omega) =\prod_{p<z,p \in\cP} \left( 1 - \frac{\omega(p)}{p}\right),
$$
we obtain
$$
|\siv - XW(z;\omega)| \leq \left| \sum_{d | P(z)} \mu(d) R_d \right|.
$$
Then, by using the triangle inequality, we obtain the worst case remainder when all the
$\mu(d)$ terms in the sum are one
\begin{equation}
\label{eq:sieve_eratosthenes_legendre}
|\siv - XW(z;\omega)| \leq \sum_{d | P(z)} |R_d|.
\end{equation}
This is the Sieve of Eratosthenes--Legendre, which crudely bounds
the error between the actual sieved set and the 
approximations generated by choosing $X$ and $\omega(p)$.
Sometimes, even approximating this sieve proves too difficult, so
the problem is weakened before searching for an upper or lower bound.
For example, see (\ref{eq:lower_bound_D(N)}).
For an upper bound, we search for functions
$\mu^+(d)$ that act as upper bounds to $\mu(d)$, in the sense that
\begin{equation}
\label{eq:mobius_upper}
\sum_{d|(n,P(z))} \mu^+(d) \geq \sum_{d|(n,P(z))} \mu(d) =
\begin{cases}
1 & (n,P(z)) = 1, \\
0 & (n,P(z)) > 1.
\end{cases}
\end{equation}
Then
\begin{equation}
\label{eq:sieve_upper_bound}
\siv \leq X \sum_{d|P(z)}\frac{\mu^+(d)\omega(d)}{d} + \sum_{d|P(z)} |\mu^+(d)||R_d|.
\end{equation}
Minimising the upper bound while ensuring $\mu^+(d)$ satisfies 
(\ref{eq:mobius_upper}) and has sufficiently small support (to reduce the size of the summation) is, in general,
very difficult.
Nevertheless, this method is used by Chen \cite{chen1973} to construct the following upper bounds of sieves 
that we have related to the Goldbach conjecture and the twin primes conjecture
\footnote{The twin primes conjecture asserts the existence of infinitely many primes $p$ for which $p+2$ is also prime.
Examples include 3 and 5, 11 and 13, 41 and 43, \ldots}.
\begin{equation}
\# \{ (p,p') : p \in \pri, p' \in \pri_{\leq 2}, p + p' = 2n \} \gg \frac{n}{\log^2 n}
\end{equation} 
\begin{equation}
\# \{ p \in \pri : p \leq x, p + 2 \in \pri_{\leq 2} \} \gg \frac{x}{\log ^2 x}
\end{equation}
where $\pri_{\leq 2}$ is the set of all integers with 1 or 2 prime factors.

The first sieve is a near miss for approximating the Goldbach sieve, 
and weakens the Goldbach conjecture to allow sums of primes, or a prime and a semiprime\footnote{A semiprime is a product of two primes.}.
The second sieve similarly weakens the twin primes conjecture. Thus, Chen's
sieves imply that there are infinitely many primes $p$ such that $p+2$ is either prime or
semiprime, and that every sufficiently large even integer $n$ can be decomposed into
the sum of a prime and a semiprime. However, the asymptotic bounds above do not show
constants that may be present, so ``sufficiently large''
could be very large indeed.  

\section{Selberg sieve}
Let $\lambda_1 = 1$ and for $d \geq 2$, let $\lambda_d \geq 2$ be arbitrary real numbers. By (\ref{eq:rewrite_sieve}), we can construct the upper bound
\begin{equation}
\cS(\cA,\cP,z) 
= \sum_{n \in \cA} \sum_{\substack{ d | n \\ d | P(z)}} \mu(d) 
\leq \sum_{n \in \cA} \left( \sum_{\substack{ d | n \\ d | P(z)}} \lambda_d \right)^2.
\end{equation}
By squaring and rearranging the order of summations \cite{halberstam1974sieve}, we obtain
\begin{equation}
\cS(\cA,\cP,z) 
\leq \sum_{\substack{d_v | P(z) \\ v = 1,2}} \lambda_{d_1} \lambda_{d_2}
	\sum_{\substack{a \in \cA \\ a \equiv 0 \mod \lcm(d_1,d_2) }} 1.
\end{equation}
This is used to construct an upper bound for the sieve in the same form as (\ref{eq:sieve_upper_bound}).
\begin{equation}
\label{eq:selberg_double_sieve}
\siv \leq X \sum_{\substack{d_v | P(z) \\ v = 1,2}} \lambda_{d_1} \lambda_{d_2} \frac{\omega(D)}{D}
 + \sum_{\substack{d_v | P(z) \\ v = 1,2}} |\lambda_{d_1} \lambda_{d_2}  R_D| = X \Sigma_1 + \Sigma_2,
\end{equation}
Now one can choose the other remaining constants $\lambda_d, d \geq 2$ so that the main term $\Sigma_1$ is as small 
an upper bound as possible, while ensuring $\Sigma_1$ dominates $\Sigma_2$.
Since this is in general very difficult even for simple sequences $\cA$, Selberg looked at the more restricted case of setting all
the constants
$$
\lambda_d = 0 \text{ for } d \geq z
$$
and then choosing the remaining $\lambda_d$ terms to minimise $\Sigma_1$, which is now a quadratic in $\lambda_d$.
Having fewer terms to deal with makes it easier to control, and hence to bound the size of the remainder term $\Sigma_2$.

  \chapter{Examining Wu's paper}\label{cha:wu_integral}

In this chapter we look at the section of Wu's paper needed to compute $C^*$.
\newcommand{\id}[2]{\1_{[#1,#2]}}
\newcommand{\ida}[2]{\1_{[\alpha_{#1},\alpha_{#2}]}}
\newcommand{\loge}[1]{\log \left( #1 \right)}
\newcommand{\ud}{\underline{d}}
\section{Chen's method}
%
Chen wished to take the set
$$
\cA = \{ N - p : p \leq N \},
$$
and apply a sieve to it, keeping only the integers in $\cA$ that are prime.
This leaves the set of all primes of the form $N-p$, which (due to double counting)
is asymptotic to $2D(N)$.

The goal is to approximate the sieve $\siv$, as described in Chapter \ref{cha:sieve}.
Chen uses the sieve Selberg \cite{selberg1952} used to prove $C^* \leq 16 + \eps$,
which has some extra constraints on the multiplicative function $\omega(p)$ in the main term, to make it
easier to estimate.
Define
\begin{equation}
V(z) = \prod_{p<z} \left( 1 - \frac{\omega(p)}{p}\right)
\end{equation}
and suppose there exists a constant $K > 1$ such that
\begin{equation}
\frac{V(z_1)}{V(z_2)} \leq \frac{\log z_2}{\log z_1} \left( 1 + \frac{K}{\log z_1}\right) \text{ for } z_2 \geq z_1 \geq 2.
\end{equation}
Then the Rosser--Iwaniec  \cite{iwaniec1980new} linear sieve is given by
\begin{align}
\label{eq:iwaniec_sieve}
\siv & \leq XV(z) \left\{ F\left( \frac{\log Q}{\log z}\right) + E \right\} + \sum_{l<L} \sum_{q | P(z)} \lambda_l^+(q) r(\cA,q), \\
\siv & \geq XV(z) \left\{ f\left( \frac{\log Q}{\log z}\right) + E \right\} - \sum_{l<L} \sum_{q | P(z)} \lambda_l^-(q) r(\cA,q),
\end{align}
where $F$ and $f$ are the solutions of the following coupled differential equations,
\footnote{The astute reader will note the similarity between $F,f$ and the Buchstab function $\omega(u)$ (see \ref{eq:buchstab})}
\begin{align}
\label{eq:F}
F(u) = 2e^{-\gamma}/u,& & & f(u) = 0, & (0 < u \leq 2)\\
\label{eq:f}
(uF(u))' = f(u-1),& & & (uf(u))' = F(u-1), & (u \geq 2)
\end{align}
The Rosser--Iwaniec sieve is a more refined version of Selbergs sieve, as the error terms $\lambda^+_l$ and $\lambda^-_l$ have some restrictions that, informally,
$\lambda^+_l$ (or $\lambda^-_l$) can be decomposed into the convolution of two other functions $\lambda = \lambda_1 * \lambda_2$.
We will be concerned only with (\ref{eq:iwaniec_sieve}), as we only need to find an upper bound for the sieve.
Lower bounds on the sieve would imply the Goldbach conjecture, which would be difficult.
Chen improved on the sieve (\ref{eq:iwaniec_sieve}) by introducing two new functions
$H(s)$ and $h(s)$ such that (\ref{eq:iwaniec_sieve}) holds with $f(s) +h(s)$ and $F(s) -  H(s)$ in place of $f(s)$ and $F(s)$ respectively \cite{wu2004chen}.
\begin{equation}
\label{eq:upper_sieve}
\siv \leq XV(z) \left\{\left( F(s) -  H(s) \right) \left( \frac{\log Q}{\log z}\right) + E \right\} + \text{error},
\end{equation}
Chen proved that $h(s) > 0$ and $H(s)> 0$ (which is obviously a required property, as otherwise these functions would make the bound on $\siv$ worse) using three set of
complicated inequalities (the largest had 43 terms!).

\section{Wu's improvement}

Wu followed the same line of reasoning as Chen, but created a new set of inequalities that describe the functions $H(s)$ and $h(s)$. Wu's inequalities 
are both simpler (having only 21 terms) and make for a tighter bound
on $\siv$. We will not look into the inequalities, but rather a general overview of the method Wu used. For full details, see \cite[p. 233]{wu2004chen}.

Let $\delta > 0$ be a sufficiently small number, $k \in \N$ and $0 \leq i \leq k$. 
Define
\begin{equation}
Q =N^{1/2 - \delta}, \quad \underline{d} =Q/d, \quad \cL = \log N, \quad  W_k = N^{d^{1+k}}.
\end{equation}
Let $\Delta \in \R$ such that $1 + \cL^{-4} \leq \Delta \leq  1 + 2 \cL^{-4}$.

Define
\begin{equation}
A(s) = sF(s)/2e^{-\gamma}, \quad a(s) = sf(s)/2e^{\gamma},
\end{equation}
where $F$ and $f$ are defined by (\ref{eq:F}) and (\ref{eq:f}).

\newcommand{\li}{\text{li}}

Define
\begin{equation}
\label{eq:sum_sieve}
\Phi(N,\sigma,s) = \sum_d \sigma(d) \cS(\cA_d, \{ p : (p,dN) = 1 \}, \underline{d}^{1/s}),
\end{equation}
where $\sigma$ is a arithmetical function
that is the Dirchlet convolution of a collection of characteristic functions
$$
\sigma = \1_{E_1} * \ldots * \1_{E_i},
$$
where 
$$
E_j = \{p : (p,N) = 1 \} \cap [V_j/\Delta,V_j),
$$
and $V_1,\ldots,V_i$ are real numbers satisfying a set of inequalities\cite[p. 244]{wu2004chen}.
Informally, the inequalities state that the $V_i$ terms are ordered in size, bounded below by $W_k$
$$
V_1 \geq \ldots \geq V_i \geq W_k,
$$
and that no one $V_i$ term can be too large. Each $V_i$ term is bounded above by both $Q$ and the previous terms $V_1,\ldots,V_{i-1}$.
$$
V_i \leq \sqrt{\frac{Q}{V_1\ldots V_{i-1}}}.
$$
So if any one $V_i$ term is big, the rest of the terms beyond will be constrained to be small.

$\Phi$ can be thought of as breaking up the sieve $\siv$ into smaller parts,
where the set that is being sieved over is only those elements in $\cA$ that are multiples of $d$,
and the index of summation is given by the $\sigma$, as $\sigma(d)$ will be zero everywhere expect on its set of support.
Breaking up the sieve this was allows Wu to prove some weighted inequalities reated to $\Phi$, that would ordinarily be too difficult to prove in general
for the entire sieve.

Define
\begin{equation}
\label{eq:thicc_theta}
\Theta(N,\sigma) = 4 \Li (N) \sum_d \frac{\sigma(d) C_{dN}}{\varphi(d) \log \underline{d}},
\end{equation}
where $C_N$ is defined in (\ref{eq:hardyconj}).

Now for $k \in \N^+, N_0 \geq 2$ and $s \in [1,10]$ we defined $H_{k,N_0}(s)$ and $h_{k,N_0}$
as the supremum of $h$ such that for all $N \geq N_0$ and functions $\sigma$ comprised of the convolution of
no more than $k$ characteristic functions, the following inequalities hold
\begin{equation}
\Phi(N,\sigma,s) \leq \{ A(s) - h\} \Theta(N,\sigma), \qquad \Phi(N, \sigma ,s) \geq \{a(s) + h \} \Theta(N, \sigma).
\end{equation}
Wu shows that both $H_{k,N_0}(s)$ and $h_{k,N_0}$ are decreasing, and defines
\begin{equation}
H(s) = \lim_{k \to \infty} \left( \lim_{N_0 \to \infty} H_{k,N_0}(s) \right),
\qquad
h(s) = \lim_{k \to \infty} \left( \lim_{N_0 \to \infty} h_{k,N_0}(s) \right)
\end{equation}
Now it is very difficult to get an explicit form of $H$ or $h$, or to even conclude anything about the behaviour beyond it is decreasing.
Wu proceeds by proving the following integral equations
\begin{equation}
\label{eq:wu_shit_int}
h(s) \geq h(s') + \int_{s-1}^{s'-1} \frac{H(t)}{t} dt, \quad H(s) \geq H(s') + \int_{s-1}^{s'-1} \frac{h(t)}{t} dt,
\end{equation}
These integral equations are still difficult to work with, as all Wu proves about $h$ and $H$ is that $H(s)$ is decreasing on $[1,10]$,
and $h(s)$ is increasing on $[1,2]$, and is decreasing on $[2,10]$.

Wu proves an upper bound for the smaller parts of the sieve $\Phi$
\begin{align*}
5\Phi(N, \sigma, s) &\leq \sum_d \sigma(d) (\Gamma_1 - \ldots - \Gamma_4 + \Gamma_5 + \ldots + \Gamma_{21}) + O_{\delta,k}(N^{1-\eta}) \\
2 \Phi(N, \sigma, s) &\leq \sum_d \sigma(d) (\Omega_1 - \Omega_2 + \Omega_3) + O_{\delta,k}(N^{1-\eta}) 
\end{align*}
where the $\Gamma_i$ terms are the dreaded 21 terms in Wu's weighted inequality \cite[p.233]{wu2004chen}. 

\section{A lower bound for $H(s)$}
Now to get an upper bound on the sieve, Wu needs to compute $H(s)$.
The integral equation is difficult to resolve, so it is weakened to give a lower
bound for $H(s)$, given the following
Wu obtains two lower bounds for $H(s)$,
\begin{proposition}
For $2 \leq s \leq 3 \leq s' \leq 5$ and $s' - s'/s \geq 2$, we have
\begin{equation}
\label{eq:wu_prop_1}
H(s) \geq \Psi_1(s) + \int_1^3 H(t) \Xi_1(t,s) \; dt,
\end{equation}

where $\Xi(t,s) = \Xi(t,s,s')$ is given by
\begin{equation}
\label{eq:Xi1}
\begin{split}
\Xi_1(t,s,s'):=
&\frac{ \sigma_0(t) }{2t} \log \left( \frac{16}{(s-1)(s'-1)} \right) 
+ \frac{ \1_{[\al_2,3]}(t)}{2t} \log \left( \frac{(t+1)^2}{(s-1)(s'-1)}\right) \\
&+ \frac{ \1_{[\al_3,\al_2]} (t)}{2t} \log \left(\frac{t+1}{(s-1)(s'-1-t)}\right),
\end{split}
\end{equation}

and where $\Psi(s)=\Psi(s,s')$ is given by

\begin{equation}
\label{eq:psi1}
\Psi_1(s,s') := \int_{2}^{s' - 1} \frac{\log(t-1)}{t} dt + 
	\frac{1}{2} \int_{1-1/s}^{1-1/s'} \frac{\log(s't-1)}{t(1-t)} dt - I_1(s,s'),
\end{equation}
where $I_1(s,s')$ is defined as
\begin{equation}
\label{eq:I1}
I_1(s,s') = \max_{\phi \geq 2} 
\iiint\displaylimits_{1/s' \leq t\leq u \leq v \leq 1/s}
\omega \left( \frac{\phi - t- u - v}{u} \right) \frac{dt \, du \, dv}{tu^2v},
\end{equation}
and $\omega(u)$ is the Buchstab Function, see (\ref{eq:buchstab}).
\end{proposition}
The equations $\Xi(s)$ and $\Psi_1(s)$ are derived from Wu's sieve inequalities, and they provide a way to rewrite the complicated
integral equations (\ref{eq:wu_shit_int}) into an inequality that can be attacked.
\begin{proposition}
For $2 \leq s \leq 3 \leq s' \leq 5$ and $s \leq k_3 \leq k_2 \leq k_1 \leq s'$ such that
\begin{align*}
s' - s'/s \geq 2,& & 1 \leq \alpha_i \leq 3 \text{ for } (1 \leq i \leq 9),& &\alpha_1 < \alpha_4,& &\alpha_5 < \alpha_8
\end{align*}
are all satisfied, then
\begin{equation}
\label{eq:wu_prop_2}
H(s) \geq \Psi_2(s) + \int_1^3 H(s) \Xi_2(t,s)\; dt,
\end{equation}
\end{proposition}

where $\Psi_2(s) = \Psi_2(s,s',k_1,k_2,k_3)$ is given by

\begin{equation}
\label{eq:psi2}
\begin{split}
\Psi_2(s,s',k_1,k_2,k_3)
=& -\frac{2}{5} \int_2^{s'-1} \frac{\log (t-1)}{t} \, dt -\frac{2}{5}  \int_2^{k_1-1} \frac{\log (t-1)}{t} \, dt \\
& - \frac{1}{5} \int_2^{k_2-1} \frac{\log(t-1)}{t} \, dt + \frac{1}{5} \int_{1-1/s}^{1-1/s'} \frac{\log(s't-1)}{t(1-t)} \, dt \\
& + \frac{1}{5} \int_{1- 1/k_3}^{1-1/k_1} \frac{\log(k_1 t - 1)}{t(1-t)} \, dt - \frac{2}{5} \sum_{i=9}^{25} I_{2,i}(s).
\end{split}
\end{equation}
and where $I_{2,i}(s) = I_{2,i}(s,s',k_1,k_2,k_3)$ is given by

\newcommand{\ds}{\mathbb{D}}
\newcommand{\intd}[1]{\int\displaylimits_{#1}}
\newcommand{\buch}[1]{\omega \left( #1\right)}
\begin{equation}
\label{eq:I2}
\begin{split}
I_{2,i}(s) &= \max_{\phi \geq 2} \intd{\ds_{2,i}} \buch{\frac{\phi-t-u-v}{u}} \frac{dt \, du \, dv}{tu^2v} \qquad \qquad \qquad (9 \leq i \leq 15), \\
I_{2,i}(s) &= \max_{\phi \geq 2} \intd{\ds_{2,i}} \buch{\frac{\phi-t-u-v-w}{v}} \frac{dt \, du \, dv \, dw}{tuv^2w} \qquad (16 \leq i \leq 19), \\
I_{2,20}(s) &= \max_{\phi \geq 2} \intd{\ds_{2,20}} \buch{\frac{\phi-t-u-v-w-x}{w}} \frac{dt \, du \, dv \, dw \, dx}{tuv^2x}, \\
I_{2,21}(s) &= \max_{\phi \geq 2} \intd{\ds_{2,21}} \buch{\frac{\phi-t-u-v-w-x-y}{x}} \frac{dt \, du \, dv \, dw \, dx \, dy}{tuvwx^2y}.
\end{split} 
\end{equation}

and the domains of integration are 
\newcommand{\D}[1]{\mathbb{D}_{2,#1}}

\begin{equation}
\begin{split}
\label{eq:integration_domains}
\D{9} &= \{(t,u,v) : 1/k_1 \leq t \leq u \leq v \leq 1/k_3 \}, \\
\D{10} &= \{(t,u,v) : 1/k_1 \leq t \leq u \leq 1/k_2 \leq v \leq 1/s \}, \\ 
\D{11} &= \{(t,u,v) : 1/k_1 \leq t \leq 1/k_2 \leq u \leq v \leq 1/k_3 \}, \\ 
\D{12} &= \{(t,u,v) : 1/s' \leq t \leq u \leq 1/k_1, 1/k_3 \leq v \leq 1/s \}, \\ 
\D{13} &= \{(t,u,v) : 1/s' \leq t \leq 1/k_1 \leq u \leq 1/k_2 \leq v \leq 1/s \}, \\ 
\D{14} &= \{(t,u,v) : 1/s' \leq t \leq 1/k_1, 1/k_2 \leq u \leq v \leq 1/s \}, \\ 
\D{15} &= \{(t,u,v) : 1/k_1 \leq  t \leq  1/k_2 \leq  u \leq  1/k_3 \leq v\leq 1/s \},\\
\D{16} &= \{(t,u,v,w) : 1/k_2 \leq t \leq u \leq v \leq w \leq 1/k_3 \}, \\ 
\D{17} &= \{(t,u,v,w) : 1/k_2 \leq t \leq u \leq v \leq 1/k_3 \leq w \leq 1/s \}, \\ 
\D{18} &= \{(t,u,v,w) : 1/k_2 \leq  t \leq  u \leq  1/k_3 \leq  v \leq  w \leq  1/s \}, \\ 
\D{19} &= \{(t,u,v,w) : 1/k_1 \leq  t \leq  1/k_2, 1/k_3 \leq u \leq  v \leq  w \leq  1/s  \}, \\ 
\D{20} &= \{(t,u,v,w,x) : 1/k_2 \leq  t \leq  1/k_3 \leq  u \leq  v \leq  w \leq  x \leq  1/s \}, \\ 
\D{21} &= \{(t,u,v,w,x,y) : 1/k_3 \leq  t \leq  u \leq  v \leq  w \leq  x \leq  y \leq 1/s \}. \\  
\end{split}
\end{equation}
This large set of integrals is derived by finding an integral equation that provides an upper bound for each term of the form
$$
\sum_d \sigma(d) \Gamma_i
$$
Taking the sum of all these integrals will give a bound for $\Phi$, and hence for $H(s)$.
The $\alpha_i$ terms are given by
\begin{align*}
\alpha_1 & := k_1 -2,  			&\alpha_2 := s'-2,\\
\alpha_3 & := s' - s'/s - 1, 	&\alpha_4 := s' - s'/k_2 - 1, \\
\alpha_5 & := s' - s'/k_3 - 1,	&\alpha_6 := s'- 2s' / k_2, \\
\alpha_7 & := s' - s'/k_1 - s'/k_3, &\alpha_8 := s' - s'/k_1 - s'/k_2 \\
\alpha_9 & := k_1 - k_1/k_2 - 1.
\end{align*}

The function $\Xi_2(t,s) = \Xi_2(t,s,s',k_1,k_2,k_3)$ is given by
\begin{equation}
\label{eq:Xi2}
\begin{split}
\Xi_2(t;s) := \frac{\sigma_0(t)}{5t} 
&\log \left( \frac{1024}{(s-1)(s'-1)(k_1 - 1)(k_2-1)(k_3-1)}\right) \\
&+ \frac{\1_{[\al_2,3]}(t)}{5t} \log \left( \frac{(t+1)^5}{(s-1)(s'-1)(k_1-1)(k_2-1)(k_3-1)}\right) \\
&+ \frac{\ida{9}{1}(t)}{5t} \loge{\frac{t+1}{(k_2-1)(k_1-1-t)}} \\
&+ \frac{\ida{5}{2}(t)}{5t} \loge{\frac{t+1}{(k_3-1)(s'-1-t)}} \\
&+ \frac{\ida{3}{2}(t)}{5t} \loge{\frac{t+1}{(s-1)(s'-1-t)}} \\
&+ \frac{\ida{1}{2}(t)}{5t} \loge{\frac{(t+1)^2}{(k_1-1)(k_2-1)}} \\
&+ \frac{\ida{7}{5}(t)}{5t(1+t/s')} \loge{\frac{{s'}^2}{(k_1 s' - s' - k_1 t)(k_3 s' - s' - k_3 t)}} \\
&+ \frac{\ida{5}{8}(t)}{5t(1-t/s')} \loge{\frac{s'(s'-1-t)}{k_1s'-s'-k_1t}} \\
&+ \frac{\ida{6}{8}(t)}{5t(1-t/s')} \loge{\frac{s'}{k_2 s' - s' - k_2t}} \\
&+ \frac{\ida{8}{2}(t)}{5t(1-t/s')} \loge{s' - 1 - t}. \\
\end{split}
\end{equation}
$\Xi_2$, together with the $\alpha$ terms, are derived by a complicated Lemma \cite[p. 248]{wu2004chen} relating three
separate integral inequality equations.

where
\begin{equation}
\label{eq:sigma_func}
\sigma(a,b,c) := \int_a^b \log \frac{c}{t-1} \frac{1}{t} \; dt,
\qquad
\sigma_0(t) := \frac{\sigma(3,t+2,t+1)}{1-\sigma(3,5,4)}.
\end{equation}

where

All of these complicated integrals are a way of bounding Wu's 21 term inequality.
For each term $\Gamma_i$ in the inequality, it has a corresponding integral.
For example, here is one term from the inequality,
\begin{equation}
\Gamma_{10} = \underset{\ud^{1/k_1} \leq p_1 \leq p_2 \leq \ud^{1/k_2} \leq p_3 \leq \ud^{1/s}}{\sum\sum\sum}
\cS(\cA_{dp_1p_2p_3},\{ p: (p,dN) = 1 \},p_2 \})
\end{equation}
which corresponds to
$$
I_{2,10}(s) = \max_{\phi \geq 2} \intd{1/k_1 \leq t \leq u \leq 1/k_2 \leq v \leq 1/s} \buch{\frac{\phi-t-u-v}{u}} \frac{dt \, du \, dv}{tu^2v}
$$
The sieve is related to the function $F$, as $F$ is part of the upper bound for the sieve (\ref{eq:upper_sieve}).
Now if we define the Dickman function $\rho(u)$ by
\begin{align}
\label{eq:dickman}
\begin{aligned}
&\rho(u) = 1 &&  1 \leq u \leq 2, \\ 
&(u-1) \rho'(u) = -\rho(u-1) && u \geq 2.
\end{aligned}
\end{align}
then we can actually write $F$ and $f$ in terms of $\omega$ and $\rho$ \cite{halberstam1974sieve}.
\begin{align}
F(u) &= e^\gamma \left( \omega(u) + \frac{\rho(u)}{u}\right) , u > 0, \\
f(u) &= e^\gamma \left( \omega(u) - \frac{\rho(u)}{u}\right) , u > 0. 
\end{align}
which provides a way to link the Buchstab function back to the sieve.
\section{Discretising the integral}
Wu comments that obtaining an exact solution is very difficult,
and provides a lower bound on $H(s)$ by splitting up the integrals in \ref{eq:wu_prop_1}
and \ref{eq:wu_prop_2} into 9 pieces.
By letting
$s_0 :=1$
and
$s_i :=2.1 + 0.1i$ for $i = 1 , \ldots, 9$,
and the fact that $H(s)$ is decreasing on the interval $[1,10]$, Wu obtains
\begin{equation}
\label{wu_prop1_num}
H(s_i)  \geq \Psi_2(s_i) + \sum_{j=1}^9 a_{i,j}  H(s_j)
\end{equation}
where
$$
a_{i,j} := \int_{s_{j-1}}^{s_j} \Xi_2(t,s_i) dt \quad i = 1,\ldots,4; \; j=1,\ldots,9
$$
and
\begin{equation}
\label{wu_prop2_num}
H(s_i)  \geq \Psi_1(s_i) + \sum_{j=1}^9 a_{i,j}  H(s_j)
\end{equation}
where
$$
a_{i,j} := \int_{s_{j-1}}^{s_j} \Xi_1(t,s_i) dt \quad i = 5,\ldots,9; \; j=1,\ldots,9
$$
Now $H(s_i)$ is given as a linear combination of all the other $H(s_j), 1 \leq j \leq 9$ values,
which simplifies the problem from resolving a complicated integral equation, to a simple linear optimisation problem.
These discretisations of the integrals can be written as matrix equations. 
\begin{equation}
\label{eq:wu_matrix_def}
\bf{A}:=\begin{bmatrix}
a_{1,1} & \ldots & a_{1,9} \\
 \vdots &\ddots & \vdots \\
 a_{9,1} & \ldots & a_{9,9} \\
\end{bmatrix},
\quad
\bf{H} := \begin{bmatrix}
H(s_1) \\
\vdots \\
H(s_9) \\
\end{bmatrix},
\quad
\bf{B} := \begin{bmatrix}
\Psi_2(s_1) \\
\vdots \\
\Psi_2(s_4) \\
\Psi_1(s_5) \\
\vdots \\
\Psi_1(s_9) \\
\end{bmatrix}
\end{equation}

Thus \ref{wu_prop1_num} and \ref{wu_prop2_num} can be rewritten as
\begin{equation}
\label{eq:wu_matrix}
\bf{(I-A)H \geq B}
\end{equation}
Wu further simplifies the problem by simply changing the inequality to an equality, and solves the 
system of simultaneous linear equations. This provides a lower bound to the linear
optimisation problem given.
So the equation below is solved for $\bf{X}$, which is easy to do.
\begin{equation}
\label{eq:wu_matrix2}
\bf{(I-A)X = B}
\end{equation}
Thus, we obtain that 
\begin{equation}
\label{eq:HgeqX}
\bf{H \geq X}
\end{equation}
Chen's constant is then equal to the largest element in the vector $\bf{X}$,
as that corresponds to a lower bound for $H(s)$ at some point $s$.
Bounds on $H$ give us bounds on the sieve, by (\ref{eq:upper_sieve}).
Naturally, we choose the best element in $\bf{X}$ (i.e.\ the largest) to get the best bound on $\bf{H}$.
It happens to be that the first element in $\bf{X}$ is always the biggest, as the function $\Psi_2(s)$ is maximised at approximately $s= 2.2$. 
Now, by taking $\sigma = \{1\}$ and $s = 2.2$ in the sum of sieves $\Phi(N,\sigma,s)$ (\ref{eq:sum_sieve})
\begin{align*}
\Phi(N,\{1\},2.2) &= \sum_d \sigma(d) \cS(\cA_d, \{ p : (p,dN) = 1 \}, \underline{d}^{1/2.2}), \\
	&=\cS(\cA, \{ p : (p,N) = 1 \}, \underline{d}^{1/2.2}),\\
	&=\cS(\cA, \{ p : (p,N) = 1 \}, N^{(1/2 - \delta)/2.2}),
\end{align*}
By letting
$$
A(s) = sF(s)/2e^\gamma
$$
and from the definition of $\Theta(N,\sigma)$,
\begin{align*}
\Theta(N,\{1\}) &= 4 \Li(N) \sum_d \frac{C_N}{\log \underline{d}} \\
				&=  \frac{4 \Li(N) C_N}{\log N^{1/2 - \delta}} 
\end{align*}
Therefore, \cite[p. 253]{wu2004chen}
\begin{align*}
\Phi(N,\{1\},2.2) &\leq \{A(2.2) - H_{k,N_0}(2.2) \} \frac{4 \Li(N) C_N}{\log N^{1/2 - \delta}}
				  &\leq 8 \{1 - x_1 \} \Theta(N) 
\end{align*}
where $x_1$ is the first element of the vector $\bf{X}$ (\ref{eq:HgeqX}).

Thus, Wu obtains an upper bound for Chen's constant.

\chapter{Approximating the Buchstab function}\label{cha:buchstab}

In this section we examine how the Buchstab Function $\omega(u)$ (which appears in
most of Wu's integrals) is computed.

\section{Background}

The Buchstab function is defined by the following delay differential
equation\footnote{The definition of $\omega(u)$ is very similar to the other difference equations $F$ and $f$ (\ref{eq:F}) and (\ref{eq:f}).}
\begin{align}
\label{eq:buchstab}
\begin{aligned}
&\omega(u) = 1/u &&  1 \leq u \leq 2, \\ 
&(u \omega(u))' = \omega(u-1) && u \geq 2.
\end{aligned}
\end{align}
\begin{figure}[!ht]
\centering
\includegraphics[width=0.7\linewidth]{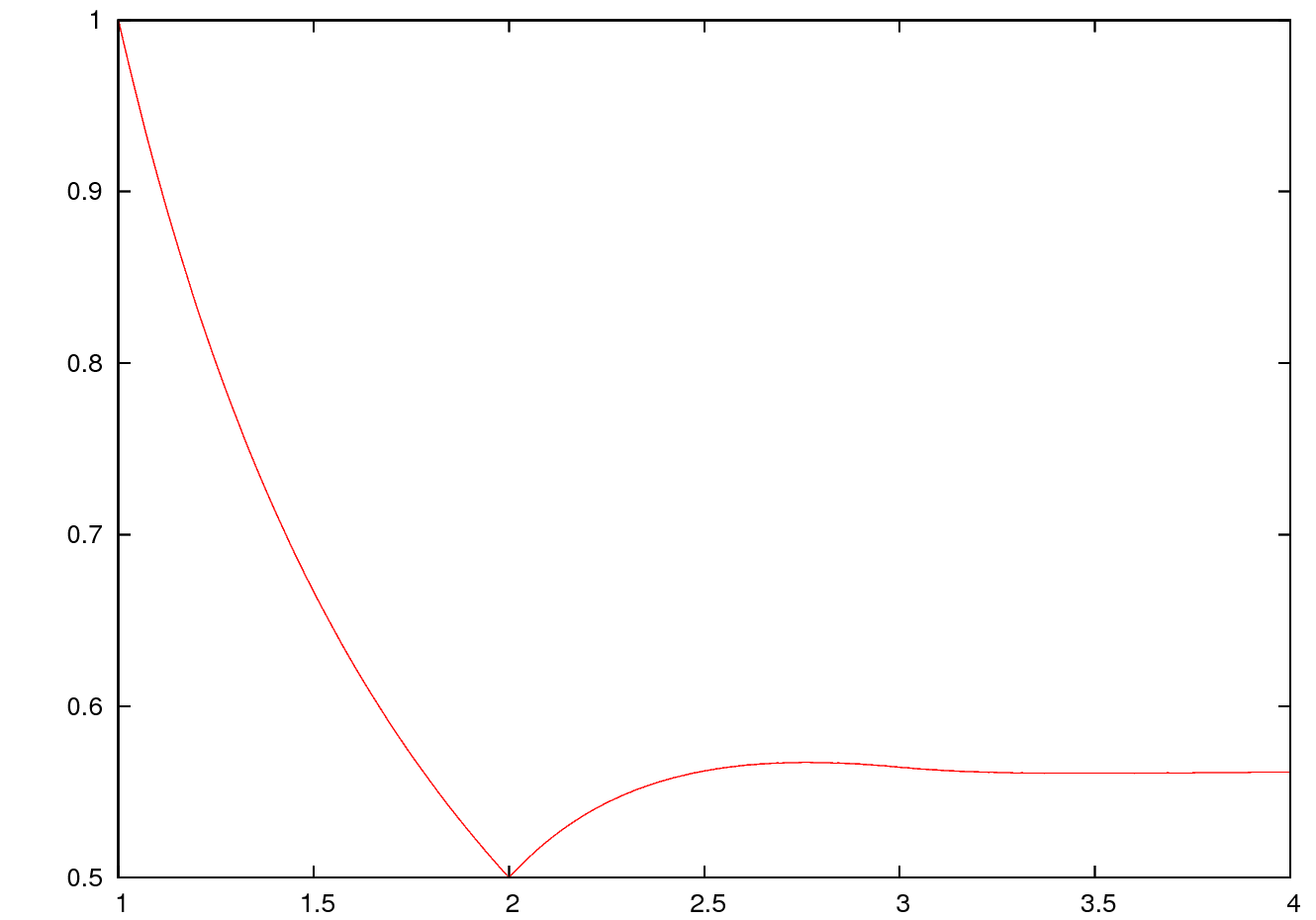}
\caption[Plot of Buchstab function]{Plot\footnotemark of the Buchstab Function $\omega(u)$ for $1 \leq u \leq 4$.}
\label{fig:Buchstab-function-graph-from-1-to-4}
\end{figure}

From the graph it appears that $\omega(u)$ quickly approaches a constant value. 
Buchstab \cite{Buc37} showed that 
\begin{equation}
\label{eq:buchstab_limit}
\lim_{u \to \infty} \omega(u) = e^{-\gamma}, 
\end{equation}
(where $\gamma$ is the Euler--Mascheroni constant)
and that the convergence is faster than exponential, i.e
\begin{equation}
\label{eq:buchstab_limit_speed}
|\omega(u) - e^{-\gamma} | = O(e^{-u}).
\end{equation}
Hua \cite{hua1951estimation} gave a much stronger bound
\begin{equation}
\label{eq:hua_bound}
|\omega(u) - e^{-\gamma} | 
\leq e^{-u(\log u + \log \log u + ( \log \log u / \log u) - 1) + O(u / \log u)}.
\end{equation}
A slightly stronger bound was obtained numerically in \emph{Mathematica} over the interval of interest,
$$
\forall u \in [2,10.5], |\omega(u) - e^{-\gamma}| \leq 0.38 e^{-1.275 u \log u}.
$$

\footnotetext{\url{https://en.wikipedia.org/wiki/Buchstab_function}}
The Buchstab function is related to rough numbers; numbers whose
prime factors all exceed some value. If we define 
\begin{equation}
\label{eq:rough_numbers}
\bold{\Phi}(x,y) = \# \{n \leq x : p|n \Rightarrow p \geq y \},
\end{equation}
i.e.\ the number of positive integers with no prime divisors below $y$, 
the following limit holds \cite{bruijn1950number}
\begin{equation}
\lim_{y \to \infty} \bold{\Phi}(y^u,y) y^{-u} \log y = \omega(u).
\end{equation}


If we define
$$
W(u) = \omega(u) - e^{-\gamma},
$$
then $W(u)$ mimics a decaying periodic function, i.e something that has similar properties to $e^{-x} \sin(x)$. The period of $W(u)$ is slowly increasing, and is expected to limit to 1.
Also, it has been shown
that in every interval $[u,u+1]$ of length 1, $W(u)$ has at least one, but no more than two, zeros. \cite{cheer1990differential} 
If we let $c_0=1$, $c_1 = 2$ and $c_2,c_3,\ldots$ denote the critical points of the Buchstab function, \cite{cheer1990differential}
we have that $c_1,c_3,c_5,\ldots$ are local minimums that are strictly increasing, and $c_2,c_4,\ldots$ are local maximums
that are strictly decreasing, in the sense that
\begin{align}
\label{eq:buchstab_critical_points}
\frac{1}{2} &= \omega(c_1) < \omega(c_3) < \omega(c_5) < \ldots \\  
1 &= \omega(c_0) > \omega(c_2) > \omega(c_4) > \ldots  
\end{align}
This allows us to obtain the trivial bound 
\begin{equation}
\label{eq:buchstab_trivial_bound}
1/2 \leq \omega(u) \leq 1.
\end{equation}
Importantly, the Buchstab function decays very quickly to a constant,
which makes approximating it easy for the purposes of computing Chen's constant.
For $1 \leq u \leq 2$, $\omega(u)$ has a closed form, and since it converges so quickly to $e^{-\gamma}$,
we can set $\omega(u)$ to be a constant when $u$ is sufficiently large. 
Cheer--Goldston published a table of the critical points and zeros of $W(u)$, and by using the critical points
$$
W(9.72844) = -9.58 \times 10^{-14}, \qquad W(10.52934) = 3.57 \times 10^{-15}
$$
we can conclude using (\ref{eq:buchstab_critical_points}) that
\begin{equation}
\label{eq:buchstab_outside}
\forall u \geq 10, |W(u)| \leq 10^{-13}.
\end{equation}
For the purposes of numerically
integrating $\omega(u)$, we only consider the region $2 \leq u \leq 10$, as beyond that, we can set 
$\omega(u) = e^{-\gamma}$ for $u \geq 10$ and bound the error by (\ref{eq:buchstab_outside}).

\section{A piecewise Taylor expansion} 
A closed form for the Buchstab function can be defined in a piecewise manner for each interval 
$n \leq u \leq n+1$, as the difference equation can be rearranged such that the value in a particular
interval depends on its predecessor. 
\begin{theorem}\cite{cheer1990differential}
If we define
$$
\omega_j(u) = \omega(u), \quad j \leq u \leq j+1
$$
and integrate (\ref{eq:buchstab}) we obtain
\begin{equation}
\label{eq:buchstab_integral}
u\omega_{j+1}(u) = \int_j^{u-1} \omega_j(t) dt + (j+1)\omega_j(j+1) \quad \text{\emph{for} } j+1 \leq u \leq j+2.
\end{equation}

\begin{proof}
$\omega(u-1)$ is defined for all $u \geq 2$, so by integrating (\ref{eq:buchstab}), 
$$
u \omega(u) = \int_2^u \omega(t-1) dt,
$$
and then by changing variables
$$
u \omega(u) = \int_1^{u-1} \omega(t) dt. 
$$
Restrict $u$ such that $j+1 \leq u \leq j+2$.
\begin{align*}
u \omega_{j+1}(u) &= \int_1^{u-1} \omega_j(t) dt \\
&= \int_j^{u-1} \omega_j(t) dt + \int_1^j \omega(t) dt \\
u \omega_{j+1}(u) - \int_j^{u-1} \omega_j(t) dt &= \int_1^j \omega(t) dt. \\
\end{align*}
Observe that $\int_1^j \omega(t) dt$ is a constant. Call this constant $K$.
$$
\forall u \in [j+1,j+2], \; K = u \omega_{j+1}(u) - \int_j^{u-1} \omega_j(t) dt.
$$
Choose $u = j+1$.
$$
K = (j+1) \omega_{j+1}(j+1).
$$
To force the Buchstab function to be continuous, we assume that the piecewise splines agree at the knots, that is,
$$
\omega_j(j+1) = \omega_{j+1}(j+1).
$$
So
$K = (j+1) \omega_{j} (j+1)$ and hence
$$
u \omega_{j+1}(u) = \int_j^{u-1} \omega_j(t) dt + (j+1)\omega_j(j+1).
$$
\end{proof}
\end{theorem}

By definition, $\omega_1(u) = u^{-1}$. To obtain $\omega_2(u)$, apply (\ref{eq:buchstab_integral}).
\begin{align*}
u \omega_2(u) &= \int_2^{u-1} \omega_1(t) dt + 2\omega_1(2) \\
			&= \int_2^{u-1} \frac{1}{t} dt + 2  \frac{1}{2}\\
			&= \log(u-1) - \log(2-1) + 1 \\
	\omega_2(u) &= \frac{\log(u-1) + 1}{u}, \quad 2 \leq u \leq 3.
\end{align*}
All values of $\omega_n(u)$ can, in principle, be obtained by repeating the above process,
but they cannot be expressed in terms of elementary functions, as
$$
\int \frac{\log(u-1) + 1}{u} du,
$$ 
is a non-elementary integral. The technique Cheer--Goldston \cite{cheer1990differential} used expresses 
each $\omega_j(u)$ as a power series about $u = j+1$.
\begin{equation}
\label{eq:buchstab_taylor}
\omega_j(u) = \sum_{k=0}^\infty a_k(j) (u - (j+1))^k.
\end{equation}
For the case of $\omega_2(u) = \frac{\log(u-1) + 1}{u}$, both $\log(u-1)$ and $u^{-1}$
are analytic, so we may write them in terms of their power series about $u=3$.
Cheer--Goldston shows that by doing this, one obtains
\begin{align}
\label{eq:buchstab_coeff_2}
a_0(2) &= \frac{1 + \log 2}{3}, \\
a_k(2) &= (-1)^{k+1} 
\left(
	-\frac{1+\log 2}{3^{k+1}} + \frac{1}{3(2^k)} \sum_{m=0}^{k-1} \frac{1}{k-m}
	\left( \frac{2}{3}\right)^m
\right).
\end{align}
We note that 
$$
\sum_{m=0}^{k-1} \frac{1}{k-m}\left( \frac{2}{3}\right)^m
\leq 
\sum_{m=0}^{k-1} \left( \frac{2}{3}\right)^m
\leq 
\sum_{m=0}^{\infty} \left( \frac{2}{3}\right)^m = 3.
$$
Whence it follows that
$$
|a_k(2)| \leq 
	\left|-\frac{1+\log 2}{3^{k+1}} + \frac{1}{2^k} \right| \leq \frac{1}{2^k}.
$$
To obtain the value of $a_k(j)$ in general, we follow the method of
Marsaglia et al.\ \cite{marsaglia1989numerical} and substitute (\ref{eq:buchstab_taylor})
into (\ref{eq:buchstab_integral}) to obtain
\begin{equation}
\label{eq:buchstab_cheer}
u \sum_{k=0}^\infty a_k(j+1) (u - (j+2))^k = \int_j^{u-1} \sum_{k=0}^\infty a_k(j) (t - (j+1))^k dt + (j+1)a_0(j). 
\end{equation}
The Taylor expansion of $\omega_2(u)$ can be shown to converge uniformly in
the interval $1.5 \leq u \leq 4.5$, as
$$
|\omega_2(u)| \leq \sum_{k=0}^{\infty} |a_k(2)| |u-3|^k
\leq \sum_{k=0}^{\infty} \frac{1}{2^k} |u-3|^k 
\leq \sum_{k=0}^{\infty} \left| \frac{u-3}{2} \right|^k
\leq \sum_{k=0}^{\infty} \left( \frac{3}{4} \right)^k = 4.
$$ 
Hence by the Weierstrass M-test\cite{goldberg1964methods}, the series for $\omega_2(u)$ is uniformly convergent for $1.5 \leq u \leq 4.5$.
This allows us to compute the integral of the Taylor series of $\omega_2(u)$, (and thereby compute $\omega_3(u)$) by swapping
the order of the sum and the integral, by the Fubini-Tonelli theorem \cite{stein2009real}. By using (\ref{eq:buchstab_integral}) we can
prove by induction that the Taylor series of $\omega_j(u)$ converges uniformly for $j \leq u \leq j+1$. This allows interchange of the sums and integrals in (\ref{eq:buchstab_cheer}). Thus, we obtain
$$
u \sum_{k=0}^\infty a_k(j+1) (u - (j+2))^k = (j+1)a_0(j) + \sum_{k=0}^\infty a_k(j)
\left( \frac{(u-(j+2))^{k+1} - (-1)^{k+1}}{k+1} \right).
$$

By making a substitution $x = u - (j+2)$ we obtain
$$
\sum_{k=0}^\infty a_k(j+1)x^{k+1} + (j+2) \sum_{k=0}^\infty a_k(j+1)x^k
= (j+1)a_0(j) + \sum_{k=1}^\infty \frac{a_{k-1}(j)}{k}\left( x^k + (-1)^{k-1} \right).
$$
By subsequently equating coefficients of like terms we obtain the following recursive 
formula for $a_k(j)$, for all $k \geq 0, j \geq 2$.
\begin{align}
\label{eq:buchstab_a}
a_0(j) &= a_0(j-1) + \frac{1}{j+1} \sum_{k=1}^\infty \frac{(-1)^k}{k+1} a_k(j-1), \\
a_k(j) &= \frac{a_{k-1}(j-1) - ka_{k-1}(j)}{k(j+1)}.
\end{align}
where the base case $a_k(2)$ and $a_0(2)$ are given in (\ref{eq:buchstab_coeff_2}).
Now one can define the Buchstab function in a piecewise manner
\begin{equation}
\omega(u) =
\begin{cases}
1/u & 1 \leq u \leq 2, \\
\omega_j(u) & j \leq u \leq j+1.
\end{cases}
\end{equation}
\subsection{Errors of Taylor expansion}
In computing $\omega(u)$, we truncate the Taylor expansion to some degree $N$, and 
then compute the coefficients of the resulting polynomials to a given accuracy \cite{cheer1990differential}.
Denote the approximation of $\omega_2(u)$ as $T_2(u)$, and define the error to be
$E_2(u) = \omega_2(u) - T_2(u)$. Define the worst case error as
$$
E = \max_{2 \leq u \leq 3} |E_2(u)|.
$$
By substuting $T_2(u)$ into
(\ref{eq:buchstab_integral}), Cheer--Goldston \cite{cheer1990differential} shows that 
one obtains a new approximation $T_3(u)$ for $\omega_3(u)$,
with a new error term of
\begin{equation}
\begin{split}
\label{eq:buchstab_error_induction}
|E_3(u)| &= \left| \frac{1}{u} \left( \int_2^{u-1} E_2(t) dt + 3E_2(3) \right)\right|\\
&\leq \left| \frac{1}{u} \left( \int_2^{u-1} E dt + 3E \right)\right|
= \frac{(u-3)E + 3E}{u} = E.
\end{split}
\end{equation}
By repeating this argument and then by induction,
$$
\forall j \geq 2 \max_{j \leq u \leq j+1} |\omega_j(u) - T_j(u)| = E.
$$
So the accuracy of $\omega_2(u)$ holds for the rest of the $\omega_k(u)$,
excluding computational errors due to machine precision arithmetic
\footnote{Real number arithmetic on a computer is implemented using floating point numbers,
which only have finite accuracy. Hence, every operation introduces rounding errors, and with enough
operations, can cause issues with the final result.}\cite{cheer1990differential}.
(If needed, \emph{Mathematica} supports arbitrary precision arithmetic, so the magnitude of the computations errors could be made arbitrarily small.)

The error $E$ is easily computed, as by the Taylor remainder theorem \cite{stewart2015calculus}, for every $u \in [2,3]$
there exists a $\xi \in [2,3]$ such that
\begin{equation}
\label{eq:buchstab_taylor_rem}
\omega_2(u) = \sum_{k=0}^N a_k(2) (u-3)^k + a_{N+1}(2)(\xi - 3)^{N+1}.
\end{equation}
Then by definition, the error term (for a Taylor expansion to order $N$) can be written as
\begin{equation}
\label{eq:buchstab_error}
E_N = \max_{2 \leq \xi \leq 3} |a_{N+1}(2) (\xi - 3)^{N+1}| \leq |a_{N+1}(2)|
\leq \frac{1}{2^{N+1}}.
\end{equation}
So if $\omega_2(u)$ is approximated with a power series up to order $N$, then the maximum
error for $\omega(u)$ anywhere is $2^{-(N+1)}$.

\subsection{Improving Cheer and Goldston's Method}

The approximation of the Buchstab function was improved, 
by computing the power series for each $\omega_j(u)$ about $j + 0.5$. This way,
$\omega_j(u)$ was evaluated at most $0.5$ away from the centre of the Taylor expansion.
In the previous case, $\omega_j(u)$ could be evaluated up to 1 away from the centre.
This provided a much lower error for the same degree Taylor expansion, improving the accuracy of
all results using the Buchstab function.
By a similar technique as above, by expanding about the middle of each interval $[j,j+1]$ we
obtain
\begin{equation}
\label{eq:buchstab_taylor_middle}
\omega_j(u) = \sum_{k=0}^\infty a_k(j) (u - (j + 1/2))^k \quad \text{for }j \leq u \leq j+1,\\
\end{equation}
where
\begin{align}
\label{eq:buchstab_middle}
a_k(2) &= (-1)^{k+1} \left(-\frac{1+\log \left(\frac{3}{2}\right)}{(5/2)^{k+1}}
+\frac{3}{5} \left( \frac{2}{3} \right)^{k+1} \sum _{m=0}^{k-1} \frac{1}{k-m} \left(\frac{3}{5}\right)^m \right),\\
a_0(j) &= \frac{1}{j+1/2} 
	\sum_{k=0}^\infty \frac{a_k(j-1)}{2^k} 
\left(
	j + \frac{(-1)^k}{2(k+1)}
\right), \\
a_k(j) &= \frac{1}{j + 1/2} \left( \frac{a_{k-1}(j-1)}{k} - a_{k-1}(j)\right).
\end{align}
Applying the Taylor remainder theorem, we calculate the error of the Taylor expansion of
$\omega_2(u)$, truncated to $N$ terms:
\begin{equation}
\label{eq:buchstab_middle_error}
E := \max_{2 \leq \xi \leq 3} |\omega_2(u) - T_2(u)| \leq |a_{N+1}(2)|| (3 - 2.5)^{N+1}|
\leq \frac{1}{2^{N+1}} |a_{N+1}(2)|.
\end{equation}
The coefficient $a_{N+1}(2)$ can be bounded above, as
$$
\sum_{m=0}^{k-1} \frac{1}{k-m} \left( \frac{3}{5} \right)^m
\leq 
\sum_{m=0}^{k-1} \left( \frac{3}{5} \right)^m
\leq
\sum_{m=0}^\infty \left( \frac{3}{5} \right)^m
=
\frac{5}{2}.
$$
So,
\begin{equation}
\begin{split}
|a_{N+1}(2)| &\leq 
\left| -\frac{1+\log \left(\frac{3}{2}\right)}{(5/2)^{N+2}}
+  \frac{3}{5} \left( \frac{2}{3} \right)^{N+2} \frac{5}{2}\right|
\leq \left( \frac{2}{3} \right)^{N+1}.
\end{split}
\end{equation}
Hence the error bound is improved:
\begin{equation}
\label{eq:buchstab_imp_middle_error}
E \leq \frac{1}{2^{N+1}} \left( \frac{2}{3} \right)^{N+1} \leq \frac{1}{3^{N+1}}.
\end{equation}
which is better than the old error bound (\ref{eq:buchstab_error}) by an exponential factor.
Again, by a similar argument to (\ref{eq:buchstab_error_induction}), this error $E$ can be shown to hold everywhere.
In practice this is a rather weak error bound, as the actual error is much less.
By computing
$$
\frac{1}{2^{N+1}} |a_{N+1}(2)|,
$$
for $10 \leq N \leq 20$, (rejecting the first few $N$ until the points settle out) and plotting the values on a log plot, we
observe the values form a line.
By fitting a curve of the form
$$
\log y = m x + c,
$$ 
to this line, the asymptotic behaviour of the numerical upper bound error is deduced to be 
approximately $O(3.31^{-N})$.
The base case for the above derivations was done with $\omega_2(u)$
instead of $\omega_1(u)$ as the Taylor expansion of $1/u$ converges slowly.
The corresponding error bounds obtained are much weaker. One could, in principle,
use $\omega_3(u)$ as the base case and thus obtain a much stronger error bound 
for $\omega_j(u), j \geq 3$,
but as mentioned above, $\omega_3(u)$ is not an elementary function.
So, $\omega_2(u)$ is the best we could hope to use. If one used the power series expression
of $\omega_3(u)$ to compute the error bounds, each $a_k(3)$ would be defined in terms of $a_k(2)$ as shown,
\begin{align*}
a_0(3) &= \frac{1}{3 + 1/2}
	\sum_{k=0}^{\infty} a_k(2) 2^{-k} 
	\left(
		j + \frac{(-1)^k}{2(k+1)}
	\right)\\
	&= 
	\frac{1}{3 + 1/2}
		\sum_{k=0}^{\infty} \frac{(-1)^{k+1}}{2^k} \left(-\frac{1+\log \left(\frac{3}{2}\right)}{(5/2)^{k+1}}
		+\frac{3}{5} \left( \frac{2}{3} \right)^{k+1} \sum _{m=0}^{k-1} \frac{1}{k-m} \left(\frac{3}{5}\right)^m \right)\\
		&\qquad \left(
			j + \frac{(-1)^k}{2(k+1)}
		\right), \\
a_k(3) &= \frac{1}{3 + 1/2} \left(\frac{a_{k-1}(2)}{k} - a_{k-1}(3)\right).
\end{align*}
and so, expanded out in full, the formulas for the coefficients would become (even more) unwieldy.
Since the error $E$ decreases exponentially quickly with respect to the degree of the polynomial,
our implementation of the Buchstab function (given above) was deemed sufficient.

  \chapter{Numerical computations}

We proceeded by implementing the above in \emph{Mathematica},
using the definition of the Buchstab function given in Chapter \ref{cha:buchstab}.
The approximation to the Buchstab function used is defined as
$$
\omega T(u) = \1_{[k+1,\infty)} e^{-\gamma} + \sum_{j=1}^k \1_{[j,j+1)} \omega_j(u),
$$
where $\omega_j(u)$ is the Taylor polynomial approximation to $\omega(u)$ around $u = j+1/2$, to some degree $N$ (\ref{eq:buchstab_taylor_middle}).
The Buchstab function $\omega(u)$ was approximated using a polynomial spline over $k$ intervals, and declared equal to the limit $e^\gamma$ (\ref{eq:buchstab_limit}) beyond the $k^{\text{th}}$ interval.

\section{Justifying Integration Method}

We attempted to compute the first of Wu's integrals (\ref{eq:I1}) without maximising with respect to $\phi$
\begin{equation}
I_1(s,s',\phi) = 
\iiint\displaylimits_{1/s' \leq t\leq u \leq v \leq 1/s}
\omega T \left( \frac{\phi - t- u - v}{u} \right) \frac{dt \, du \, dv}{tu^2v},
\end{equation}
but ran into several problems. 
Using \emph{Mathematica}'s inbuilt integration routine, the computation
would either finish quickly, or never halt, depending on the value of $\phi$ chosen.
It was discovered that the function could only be integrated for large $\phi$, such that
$\omega T(u)$ would always be in the constant region. These problems were not resolved. 
Instead, we integrated the Buchstab function by computing an anti-derivative of $\omega T(u)$,
and applying the fundamental theorem of calculus (FTC), three times.

We can show that FTC validly applies in this instance.
As $\omega T(u)$ is a piecewise spline of polynomials, it is continuous everywhere except at the points 
$u = 2, 3, \ldots, k+1$ where the splines meet. Using two theorems about Lebesgue integration \cite{stein2009real}
\begin{theorem}
If $g$ is integrable, and $0 \leq f \leq g$, then $f$ is integrable.
\end{theorem}

\begin{theorem}
If $f$ is integrable on $[a,b]$, then there exists an absolutely continuous function $F$ such that $F'(x) = f(x)$
almost everywhere, and in fact we may take $F(x) = \int_a^x f(y) \; dy$.
\end{theorem}
The Buchstab function is strictly positive and bounded, as $1/2 \leq \omega(u) \leq 1$ (\ref{eq:buchstab_trivial_bound}).
The Taylor approximation $\omega T(u)$ will also be non-negative and bounded, as we can easily bound the difference between $\omega T(u)$ and $\omega(u)$ 
to be less than any $\eps > 0$ (\ref{eq:buchstab_imp_middle_error}).
Thus if $\omega T(u)$ is defined using $k$ many splines, we can bound $\omega T(u)$ above by a constant on the interval $[1,k+1]$. Beyond $u=k+1$,
$\omega T(u)$ is constant, so $\omega T(u)$ is integrable on any interval of the form $[1,k+1+r]$ for some $r > 0$.
We can thus obtain FTC as
$$
F(x) = \int_a^x f(y) \; dy \Rightarrow F(b) - F(a) = \int_a^b f(x) \; dx.
$$
We remove the maximisation over $\phi$, and set $I_1$ to be a function dependant on $\phi$.
We can convert the domain of integration 
$$
\{ (t,u,v) : 1/s' \leq t\leq u \leq v \leq 1/s \}
$$
into three definite integrals, which is more suitable to apply FTC.
\begin{equation}
\tilde{I}_1(s,s',\phi) = \int_{v=1/s'}^{v=1/s} \int_{t=1/s'}^{t=v} \int_{u=t}^{u=v} \omega T \left(\frac{\phi-t-u-v}{u} \right) \frac{1}{t u^2 v} \; du \; dt \; dv.
\end{equation}
it was expected that the FTC could be used 3 times to compute $I_1$, however 
the values computed were nonsense. When attempting to replicate the entry in Wu's table (\ref{tab:wu}) for $i=6$,
we expected $0.0094\ldots$ and obtained a large negative result, $-10000$ or so.
It is believed that since $\omega T(u)$
was defined in a piecewise manner, integrating $\omega T \left( \frac{\phi-t-u-v}{u} \right) $ would involve resolving
many inequalities, which would only increase in complexity after three integrations. 

Attempts to numerically integrate $\omega T(u)$ were successful, however it was much faster to use \emph{Mathematica} to numerically solve the
differential delay equation defining the Buchstab function (\ref{eq:buchstab}) with standard numerical ODE solver routines, 
obtain a numerical solution for $\omega(u)$, and
numerically integrate the result. The resulting values were also closer to
Wu's. The following table includes Wu's results, compared to the results we obtained.

\newlength{\offsetpage}
\setlength{\offsetpage}{5.0cm}

\newenvironment{widepage}{\begin{adjustwidth}{-\offsetpage}{-\offsetpage}%
    \addtolength{\textwidth}{3\offsetpage}}%
{\end{adjustwidth}}

\begin{table}[ht!]
	\begin{adjustwidth}{-2em}{-2em}
		{\small
	    \begin{tabular}{ | c | c  | c | c | c | c | c | c | c | c |}
	    \hline
	    $i$ & $s_i$ & $s'_i$ & $k_{1,i}$ & $k_{2,i}$ & $k_{3,i}$ & $\Psi_1(s_i)$ (Wu)& $\Psi_1(s_i)$ & $\Psi_2(s_i)$ (Wu) & $\Psi_2(s_i)$  \\ \hline
	    1   & 2.2   & 4.54   & 3.53      & 2.90      & 2.44      &               &                 & 0.01582635     & 0.01615180   \\
	    2   & 2.3   & 4.50   & 3.54      & 2.88      & 2.43      &               &                 & 0.01224797     & 0.01547663   \\
	    3   & 2.4   & 4.46   & 3.57      & 2.87      & 2.40      &               &                 & 0.01389875     & 0.01406834   \\
	    4   & 2.5   & 4.12   & 3.56      & 2.91      & 2.50      &               &                 & 0.01177605     & 0.01187935 \\
	    5   & 2.6   & 3.58   &           &           &           & 0.00940521   & 0.00947409     &                 &                  \\
	    6   & 2.7   & 3.47   &           &           &           & 0.00655895   & 0.00659089     &                 &                  \\
	    7   & 2.8   & 3.34   &           &           &           & 0.00353675   & 0.00354796     &                 &                  \\
	    8   & 2.9   & 3.19   &           &           &           & 0.00105665   & 0.00105838     &                 &                  \\
	    9   & 3.0   & 3.00   &           &           &           & 0.00000000   & 0.00000000     &                 &                 \\
	    \hline
	    \end{tabular}
	    }
	    \caption{Table of values for $\Psi_1$ and $\Psi_2$, compared with Wu's results \cite{wu2004chen}.}\label{tab:wu}
	\end{adjustwidth}
\end{table} 

Given an $s_i$, Wu chose the parameters $s'_i,k_{1,i},k_{2,i},k_{3,i}$ to maximise $\Psi_1(s_i)$ or $\Psi_2(s_i)$.
All the values we calculated were an overshoot of Wu's results, so it was not surprising that we obtained a smaller value for Chen's constant (\ref{eq:wu_matrix}).
For $i=5, 6, \ldots,9$, the values of $s'_i$ Wu gave were verified to maximise our version of $\Psi_1(s_i)$, 
by trying all values $s'_i = 3, 3.01, \ldots, 5.00$, which is given by the constraint $3 \leq s'_i \leq 5$ on (\ref{eq:psi1}). So even though our $\Psi_1$ did not match Wu's, both
were maximised at the same points.
This indicated that we could use our ``poor mans'' $\Phi_1$ to investigate the behaviour of $\Phi_1$, but not to compute exact values for it. 

When computing the values for this table, we first wanted to see for a fixed $s$ and $s'$, how $I_1(s)$ varied as a function of $\phi$.
It was determined that the function did not grow too quickly near the maximum, and as $\phi$
grew large, $\Psi_1$ tended to a constant, which is consistent with $\omega(u)$ quickly tending to $e^{-\gamma}$.
Thus, we can apply simple numerical maximisation techniques, without  
getting trapped in local maxima, or missing the maximum because the function changes too quickly.
We only need to maximise $\phi$ over a small interval $\phi \in [2,4]$, and we
obtain the maximal value by computing $I_1$ for a few points in the interval of interest,
choosing the maximum point, and then computing again for a collection of points in a small neighbourhood about the previous maximum.
This is iterated a few times until the required accuracy is obtained. (Listing \ref{lst:I1}).
\begin{figure}
\centering
\includegraphics[width=0.7\linewidth]{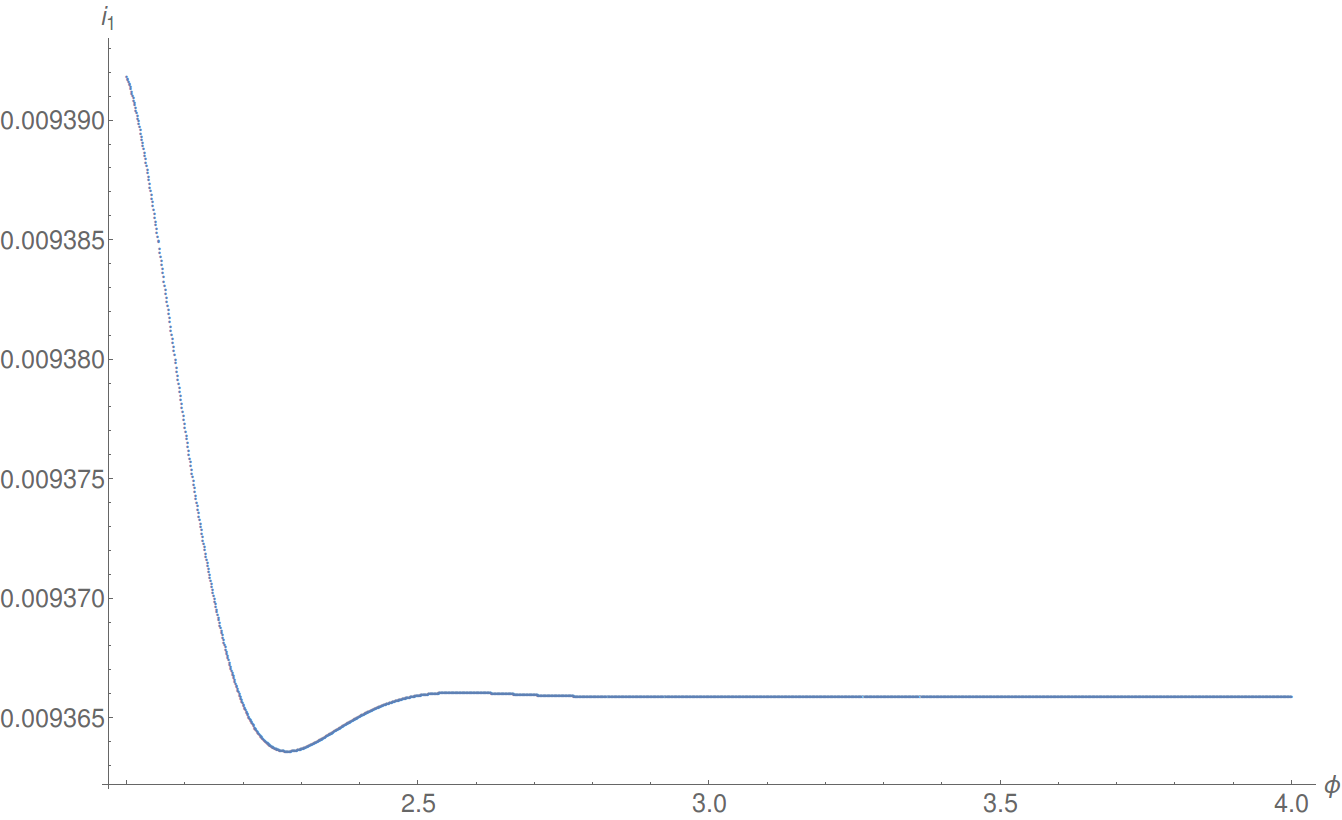}
\caption[]{Plot of $I_1(\phi,2.6,3.58)$ for $2 \leq \phi \leq 4$}
\label{tab:i1_max}
\end{figure}

\begin{minipage}{\linewidth}
\begin{lstlisting}[label={lst:I1},caption={Algorithm to compute $I_1(s,s')$},captionpos=b]
def $I_1(s,s')$:
   $\phi_{low} := 2$
   $\phi_{high} := 4$
   $\eps := 0.1$
   $\phi_{MAX} := 2$
   while($\eps \geq 0.001$):
      $\Phi := \{\phi_{low}, \phi_{low} + \eps, \phi_{low} + 2\eps, \ldots, \phi_{high}\}$
      $\phi_{MAX} := \phi \in \Phi \text{ such that } \tilde{I}(\phi,s,s') \text{ is maximal}$
	  $\phi_{low} := \max(\phi_{MAX} - \eps, 2)$
	  $\phi_{high} := \phi_{MAX} + \eps$
	  $\eps := \eps/10$
   return $\tilde{I}_1(\phi_{MAX},s,s')$
enddef
\end{lstlisting}
\end{minipage}

\begin{figure}
\centering
\includegraphics[width=0.7\linewidth]{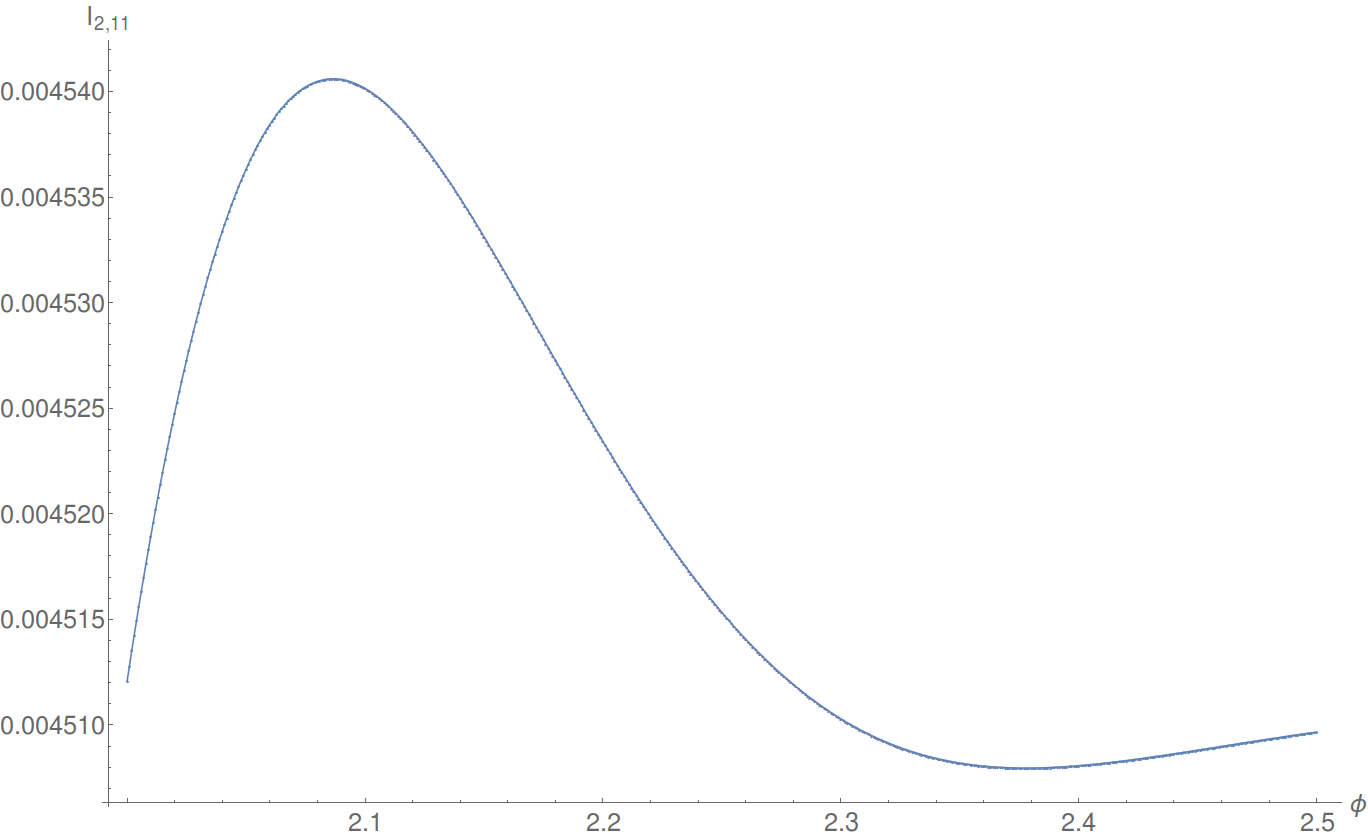}
\caption{Plot of $I_{2,11}$ for $s=2.6,s'=3.58,2\leq \phi \leq 2.5$}
\label{fig:i2-11}
\end{figure}

Initially, it seemed odd that for all the values of $s$ and $s'$ that were tried, $\tilde{I}_1(\phi,s,s')$ appeared
to always be maximal when $\phi=2$. The actual maximum for $\tilde{I}_1(\phi,s,s')$ occurs at a point $\phi < 2$,
but was being cut off by the constraint that $\phi \geq 2$. It seemed unusual to us that Wu would maximise an integral over $\phi \geq 2$,
when the integral is always maximal when $\phi = 2$. 
However, for some of the other integrals in (\ref{eq:I2}), the function was maximal at some point $\phi > 2$.
This behaviour is made clear when we plot $I_{2,11}(\phi,s=2.6,s'=3.58)$ (Figure \ref{fig:i2-11}) as a function of $\phi$, for a fixed value of $s$ and $s'$. 

For some of the integrals needed to compute $I_2$, the numerical integration did not converge. Whether
this was a property of the integrals (i.e.\ they were not defined for those particular values of $\phi$)
or a symptom of how the numerical integration routines in \emph{Mathematica} operate was not determined. To maximise the
integrals with respect to $\phi$, we could no longer assume they were defined for all $\phi \geq 2$.
Instead, it was found that for each integral, there would be a corresponding value $\phi_{low}$ such that
the method of numerical integration converged for all $\phi \geq \phi_{low}$.
This value of $\phi_{low}$ was computed by method of bisection (Listing \ref{lst:bisect}).

\begin{minipage}{\linewidth}
\begin{lstlisting}[label={lst:bisect},caption={Algorithm to compute the least $\phi$ such that $I_{2,i}(\phi,s)$ is defined},captionpos=b]
def findphi($i,s,s',k_1,k_2,k_3$):
   if(computing $I_{2,i}(2,s,s',k_1,k_2,k_3)$ suceeded):
      return $2$
   $\phi_{low} := 2$
   $\phi_{high} := 5$
   $\eps := 0.001$
   while($|\phi_{high} - \phi_{low}| \geq \eps$):
      $\phi_{mid} := (\phi_{high} + \phi_{low})/2$
      if(computing $I_{2,i}(\phi_{mid},s,s',k_1,k_2,k_3)$ suceeded):
         $\phi_{high} := \phi_{mid}$
      else:
         $\phi_{low} := \phi_{mid}$
   return $\phi_{high}$
enddef
\end{lstlisting}
\end{minipage}

$I_1$ and $I_{2,i}$ can now be computed, and thus $\Psi_1$ (\ref{eq:psi1}) and $\Psi_2$ (\ref{eq:psi2}) can
be easily computed using standard numerical integration techniques. This gives the entries of the vector $\bf{B}$ 
(\ref{eq:wu_matrix_def}).
Computing the matrix $\bf{A}$ (\ref{eq:wu_matrix_def}) is comparatively easy, as each element $a_{i,j}$ of the matrix $\bf{A}$ is the integral
of either $\Xi_1$ (see \ref{eq:Xi1}) or $\Xi_2$ (see \ref{eq:Xi1}) over an interval, which is easily computed numerically.
Thus, the system of linear equations
$
\bf{(I-A)X=B}
$
(where $\bf{I}$ is the identity matrix) can be solved for $\bf{X}$.
This is compared to the value for $\bf{X}$ that Wu obtained.
\begin{equation}
\label{eq:X_compare}
\bf{X}_{Wu} =
\left(\begin{matrix}
0.0223939\ldots\\ 
0.0217196\ldots\\ 
0.0202876\ldots\\ 
0.0181433\ldots\\ 
0.0158644\ldots\\ 
0.0129923\ldots\\ 
0.0100686\ldots\\ 
0.0078162\ldots\\
0.0072943\ldots \\ 
\end{matrix} \right)
\quad
\bf{X} =
\left(\begin{matrix}
\textbf{0.022}7656\ldots\\ 
\textbf{0.021}9942\ldots\\ 
\textbf{0.020}5028\ldots\\ 
\textbf{0.018}2930\ldots\\ 
\textbf{0.015}2937\ldots\\ 
\textbf{0.012}6404\ldots\\ 
\textbf{0.0}099076\ldots\\ 
\textbf{0.0078}020\ldots\\
\textbf{0.007}3089\ldots \\ 
\end{matrix} \right)
\end{equation}
So from these vectors, we use (\ref{eq:HgeqX}) to obtain two different values for $\K$,
ours and Wu's,
\begin{equation}
\begin{split}
\K_{Wu} &= 8(1-0.0223938) \leq 7.82085, \\
\K &= 8(1-0.0227655) \leq 7.8178752,
\end{split}
\end{equation}
which seems to be an improvement on Chen's constant, but since the 
calculations used to obtain this value are numerical in nature and don't match Wu's results, it can
be hard to conclude anything useful from them.
Nonetheless, the effect on Chen's constant of sampling $s_i$ with finer granularity was investigated.
We could no longer use the table Wu had given (\ref{tab:wu}), so we had to numerically optimise over both $\phi$ and $s'_i$.
Each evaluation of $\Psi_2$ is very expensive, as we have to compute each of the integrals given in (\ref{eq:I2}).
This means it is computationally difficult to compute $\Psi_2(s_i)$ without the corresponding parameters $\phi,s'_i,k_{1,i},k_{2,i},k_{3,i}$,
as for an interval size of $n$ points, we now need to compute $n^5$ points, which can quickly become infeasible.
We focused on increasing the sample size for $\Psi_1$, by diving the interval of $[2,3]$ into $[2,2.6] \cup [2.6,3.0]$,
and increasing the resolution of sampling in the interval $[2.6,3.0]$.
By defining
$$
s'_i = \begin{cases}
2.1 + 0.1i & 1 \leq i \leq 4, \\
2.6 + 0.01(i-5) & 5 \leq i \leq 45.
\end{cases}
$$ 
we use the same values before (See \ref{tab:wu}) for $1 \leq i \leq 4$ and then compute new values of $\Psi_1$ for each new value of $s_i$,
so we evaluate $\Psi_1$ at 40 points instead of 5.
We can then compute $\Psi_1$ by numerically maximising over $s'_i$, in a similar fashion to computing $I_1$ (See Listing \ref{lst:I1}).
Note that $\Psi_1$ depends on $I_1$, so we are optimising with respect to two variables, $\phi$ and $s'_i$.
We also need to recompute the matrix $\bf{A}$ given by (\ref{eq:wu_matrix_def}) on the new set of intervals. Thus, the component $a_{i,j}$ of $\bf{A}$ is
now given by
\begin{equation}
\label{eq:aij}
a_{i,j} :=
\begin{dcases}
\int_{s_{j-1}}^{s_j} \Xi_2(t,s_i) dt & i = 1,\ldots,4; \; j=1,\ldots,45 \\
\int_{s_{j-1}}^{s_j} \Xi_1(t,s_i) dt & i = 5,\ldots,45; \; j=1,\ldots,45
\end{dcases} 
\end{equation}
So we now have a new matrix $\bf{A}$ and vector $\bf{B}$ with which to solve $\bf{(I-A)X =B}$.
$$
\bf{A}:=\begin{bmatrix}
a_{1,1} & \ldots & a_{1,45} \\
 \vdots &\ddots & \vdots \\
 a_{45,1} & \ldots & a_{45,45} \\
\end{bmatrix}
,
\bf{B} := \begin{bmatrix}
\Psi_2(s_1) \\
\vdots \\
\Psi_2(s_4) \\
\Psi_1(s_5) \\
\vdots \\
\Psi_1(s_{45}) \\
\end{bmatrix}
\Rightarrow
\bf{X} =
\left(\begin{matrix}
0.0228801\ldots\\ 
0.0221067\ldots\\ 
0.0206136\ldots\\ 
0.0184024\ldots\\ 
0.0184024\ldots\\ 
0.0153999\ldots\\ 
0.0157765\ldots\\ 
0.0163132\ldots\\
0.0170428\ldots \\
\vdots 
\end{matrix} \right)
$$

Once that was complete, we redid the calculations again with even finer granularity for $s'_i$, choosing
$$
s'_i = \begin{cases}
2.1 + 0.1i & 1 \leq i \leq 5, \\
2.6 + 0.001(i-5) & 6 \leq i \leq 405.
\end{cases}
$$

These calculations were run in parallel on a computer equipped with an i5-3570K CPU overclocked to 3.80Ghz, and with 6GB of RAM.
\begin{table}[ht!]
\centering
\begin{tabular}{|c|c|c|c|c|c|}
\hline
Points & Time $\bf{B}$ & Time $\bf{A}$ & $x_1$ & $C^*$ \\ \hline
 Wu & - & - & 0.0223938 & 7.82085 \\
 5 & 1m27s & $<1$s & 0.0227655 & 7.8178752 \\
 40 & 11m31s & 10s & 0.0228800 & 7.81696  \\
 400 & 1h53m & 11m12s & 0.0229275 & 7.81658  \\
 \hline
\end{tabular}
\caption{Improvements on $C^*$ by diving $[2.6,3.0]$ into more points\label{tab:compute}}
\end{table}
The time taken to compute additional intervals increases rapidly,
the projected time for 4000 intervals would take longer than a day to compute, with likely minimal improvement
on the constant.

\section{Interpolating Wu's data}

Since $\Psi_2$ is complicated to evaluate, we looked at approximating it by interpolation, so see what effects
additional values of $\Psi_2$ would have on Chen's constant. We assume the best possible case where possible, to obtain the best possible value for $C^*$.
By demonstrating that even under the best assumptions the improvements on $C^*$ are minimal, this should indicate that improving $C^*$ using Wu's integrals will
give a similarly small improvement.  
$\Psi_2$ appears to be decreasing by Chen's table, and so too for $\Psi_1$, so we interpolate the points with \emph{Mathematica} in-built interpolation routine.
By plotting the interpolant of $\Psi_2$ and $\Psi_1$, we can see one overtakes the other at around $s = 2.55$. This is why Wu had split up his table in this way,
choosing $\Psi_1$ for the lower half and $\Psi_2$ for the upper, to get the best possible bound.

We define
\begin{equation}
\Psi_B = 
\begin{cases}
\Psi_2(s) & s \leq m \\
\Psi_1(s) & s \geq m
\end{cases}
\end{equation}
where $m$ is the point where $\Psi_2(m) = \Psi_1(m)$.
Thus, $\Psi_B$ is the best case scenario, by essentially taking the maximum of the two upper bounds. 
We need the values for $s',k_1,k_2,k_3$ to compute $C^*$, as we need to compute the matrix $\bf{A}$, which in turn depends on the values of $a_{i,j}$ (see \ref{eq:aij}).
So we obtain those variables by interpolation also.

We then discretise the interval $[2,3]$ from coarse to fine, and calculate the new value of $C^*$ for each of these discretisations.
In doing so, we can get an estimate as to how much additional intervals of integration impact the value of $C^*$.
After computing $C^*$ for a varying number of intervals, from Wu's original 9 up to 600 (which took 84 minutes, 40 seconds on the same hardware used to compute Table \ref{tab:compute}),
we can clearly see that adding more intervals does make a slight difference, but there are diminishing returns beyond a hundred or so. One would expect not to make any gains on $C^*$ by attempting to
compute a million intervals.

\begin{figure}[ht] 
\newcommand{\quadsize}{0.48}
  \begin{subfigure}[t]{\quadsize\linewidth}
    \centering
    \includegraphics[width=\linewidth]{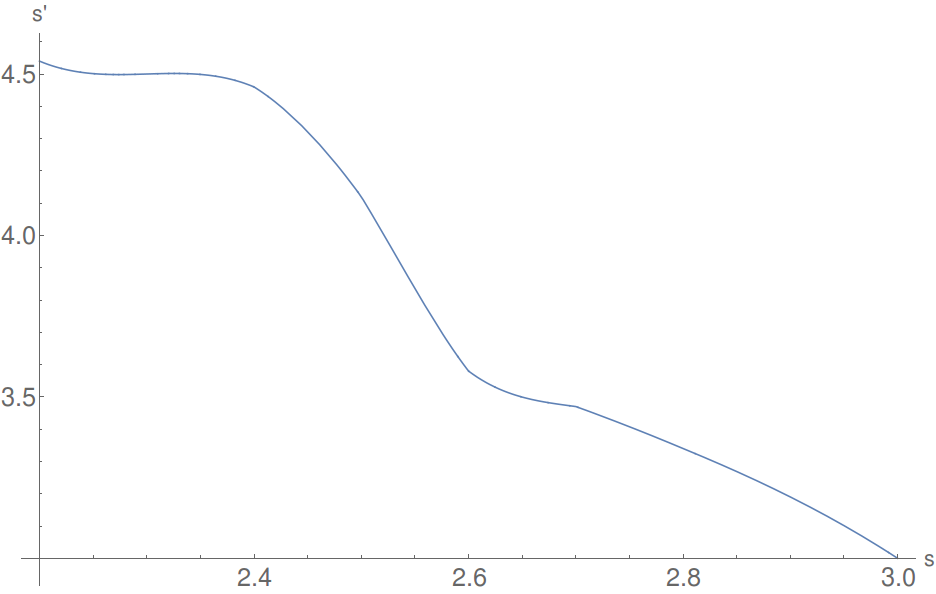} 
       \caption{Interpolant for $s'$}
       \label{fig:sint} 
    \vspace{4ex}
  \end{subfigure}
 	\hfill
  \begin{subfigure}[t]{\quadsize\linewidth}
    \centering
    \includegraphics[width=\linewidth]{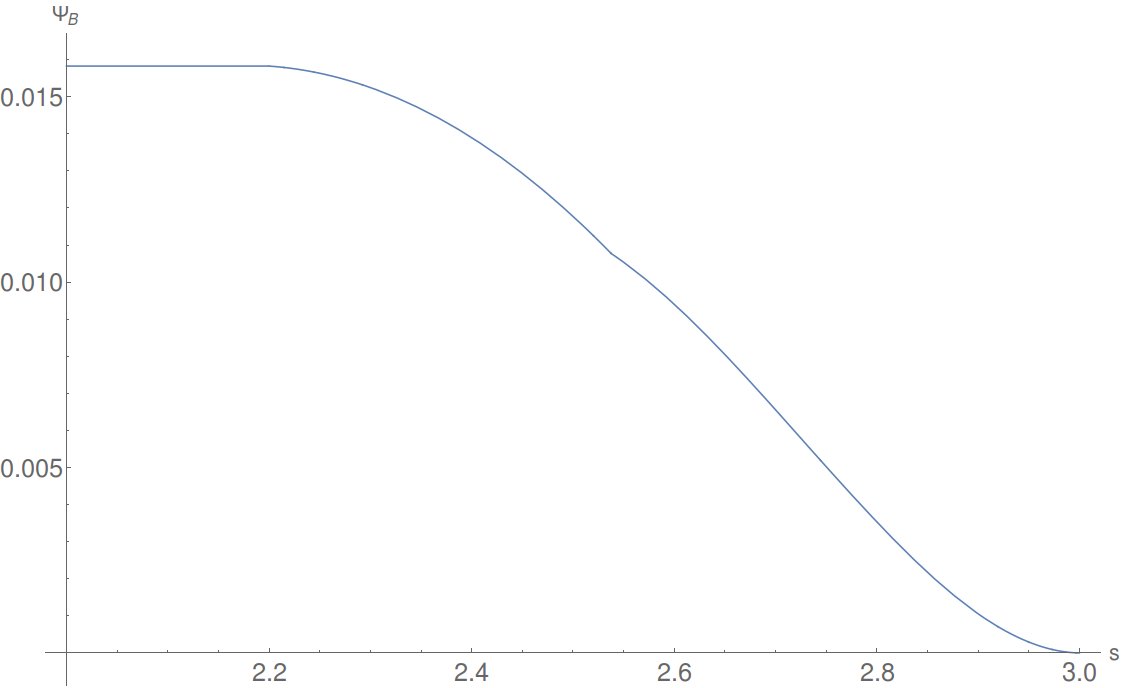} 
      \caption{Interpolant of $\Psi_B$, the best of $\Psi_1$ and $\Psi_2$ put together} 
      \label{fig:psiB}
    \vspace{4ex}
  \end{subfigure} 
  
  \medskip
  
  \begin{subfigure}[t]{\quadsize\linewidth}
    \centering
    \includegraphics[width=\linewidth]{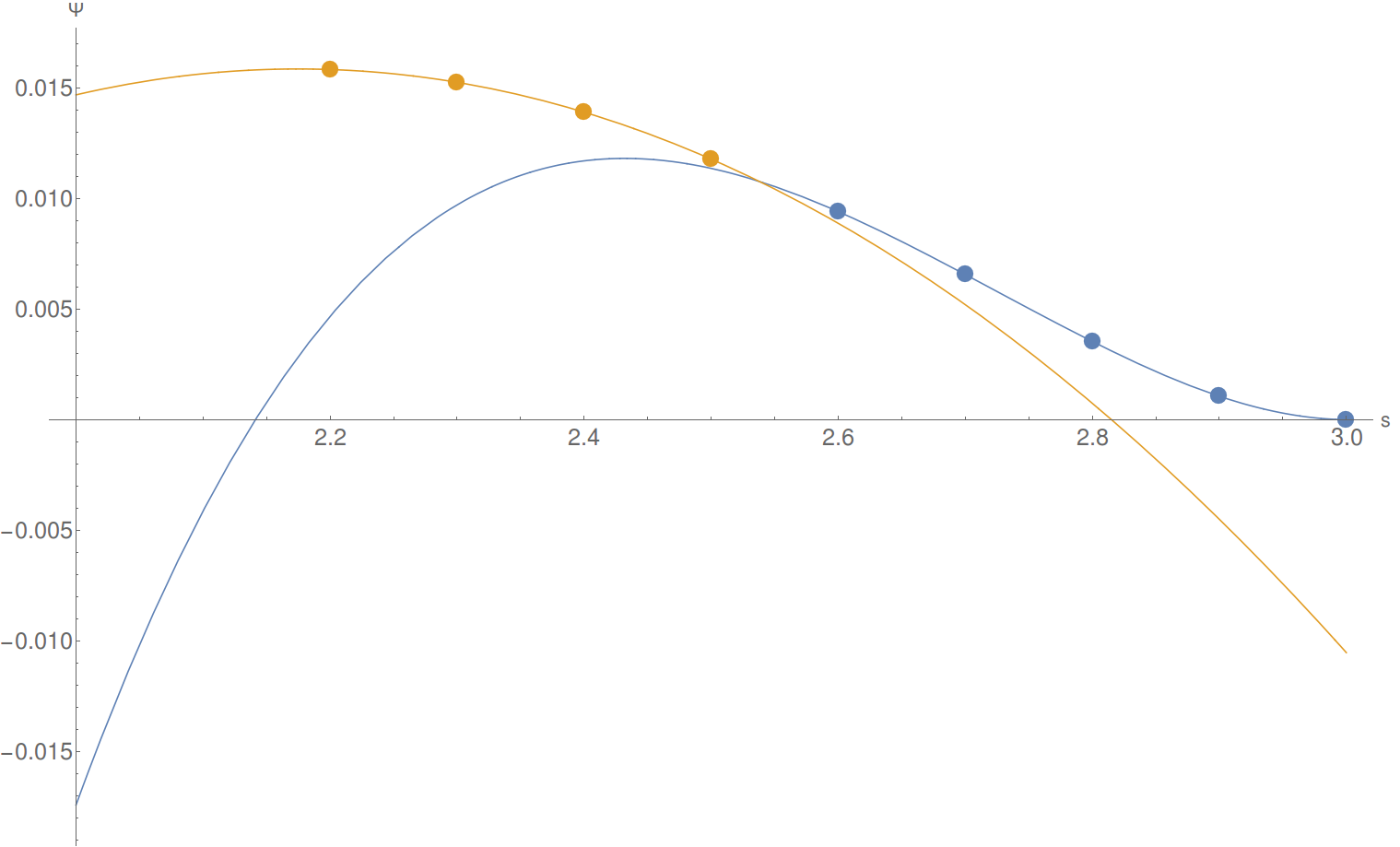} 
    \caption{Interpolant of $\Psi_2$ (Orange) and $\Psi_1$ (Blue) from Wu's data}
    \label{fig:psi_int5}
  \end{subfigure}
  \hfill
  \begin{subfigure}[t]{\quadsize\linewidth}
    \centering
	\includegraphics[width=\linewidth]{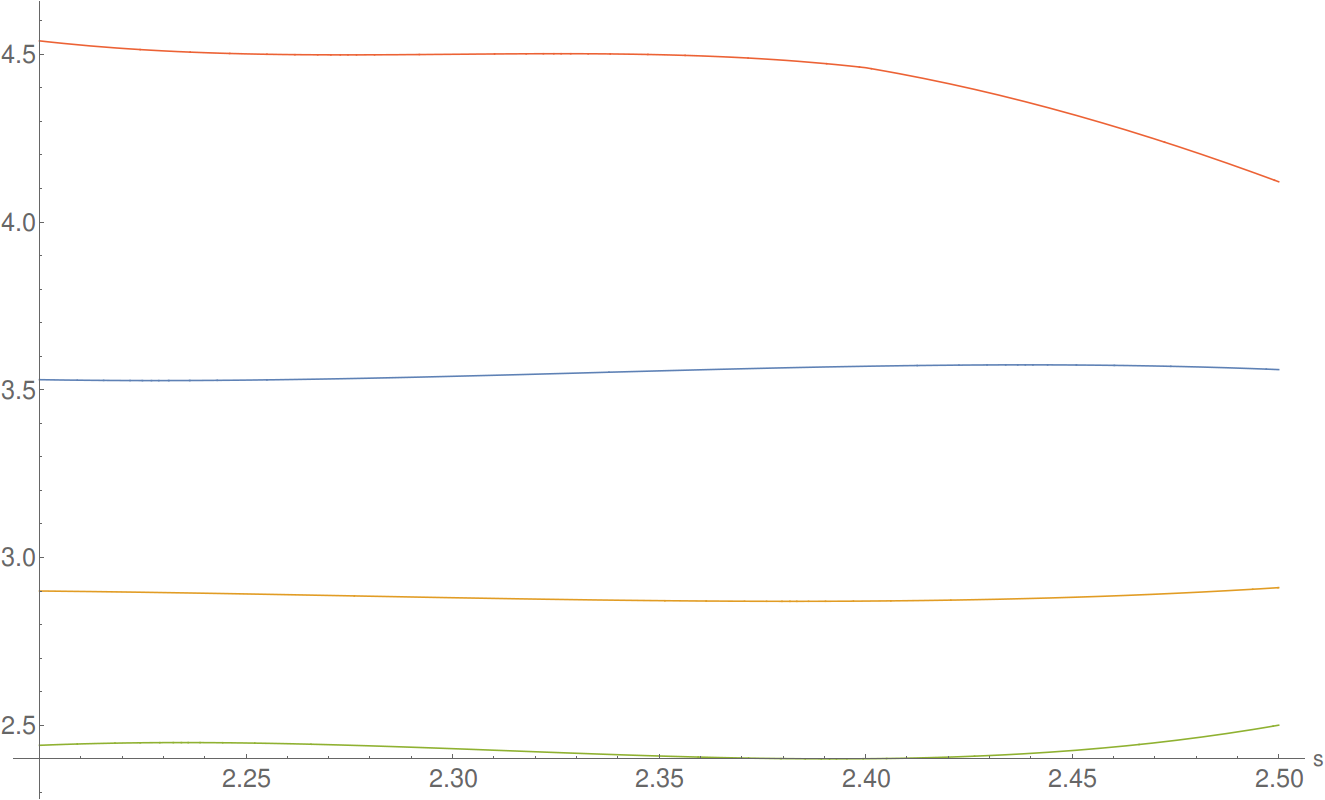} 
    \caption{Interpolants for $s',k_1,k_2,k_3$ (top to bottom)} 
  \end{subfigure} 
  \caption{Various interpolants of Wu's data}
  \label{fig:interpolant}
  \centering
       \vspace*{\floatsep}
  \includegraphics[width=0.7\linewidth]{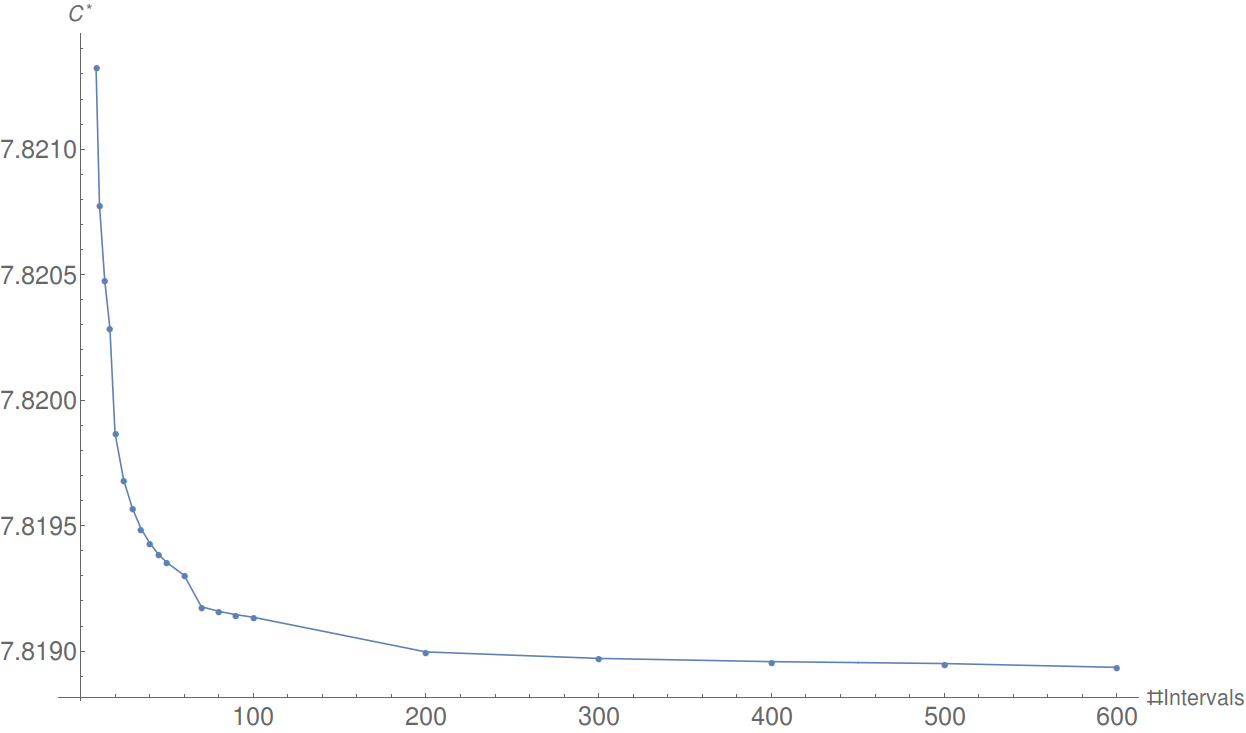}
  \caption{Value of $C^*$ versus number of intervals used in integration.}
  \label{fig:chens_int} 
\end{figure}

  \chapter{Summary}

\section{Buchstab function}

We have improved on the approximation of the Buchstab function $\omega(u)$ that Cheer--Goldston \cite{cheer1990differential} gave.
We used a spline, where each polynomial in the spline is the Taylor expansion of $\omega(u)$ 
about the centre $u = j + 0.5$ of the interval $[j,j+1]$, whereas Cheer--Golston expanded about the right endpoint $u = j+1$ of this interval. 
We derived an implicit formula
for the coefficients $a_j(k)$ for the power series
$$
\omega_j(u) = \sum_{k=0}^\infty a_j(k) (u - (j + 0.5))^k,
$$
where each $a_j(k)$ is defined in terms of the previous values $a_n(m)$ for  $2 \leq n < j, 0 \leq m < k$,
and the values of $a_2(k)$ are derived from the Taylor expansion of 
$$
\omega_2(u) = \frac{\log(u-1) + 1}{u}.
$$
The coefficients are more complicated than the ones Cheer--Goldston obtained, but
when the Taylor approximation of $\omega(u)$ is truncated to degree $N$, we show that the corresponding error is
less than $1/3^{N+1}$ (\ref{eq:buchstab_imp_middle_error}) versus Cheer--Goldston's error bound of $1/2^{N+1}$. 

Although the spline is not continuous at the knots, the size of the discontinuities can be made arbitrarily small
by increasing the degree of the power series.

\section{Chen's constant}

An attempt was made to use our approximation to the Buchstab function to evaluate the integrals Wu gave, to place an
upper bound on $C^*$. This was not successful, so instead we numerically solved the ODE that defined $\omega(u)$,
and numerically evaluated Wu's integrals instead. This gave us approximately the same results as Wu.
Wu resolved a set of functional inequalities using a discretisation of $[2,3]$ into 9 pieces, so we replicated the results, and then
increased the resolution for the section $[2.6,3.0]$ of the interval from 5 points to 40, and then 400.
The corresponding change in $C^*$ was as expected, minimal. The relative difference in $C^*$ between sampling 5 points and 400 points in the 
interval $[2.6,3.0]$ was only $1.66 \times 10^{-2} \%$. This leads us to believe that we should expect similar results, had we computed an exact form
of Wu's integrals.

We interpolated the values of $\Psi_1,\Psi_2$ and the variables required from Wu's data,
and (under the assumption that these interpolations are a reasonable approximation to the exact form),
we have shown that Wu was indeed right in asserting that cutting $[2,3]$ into more subintervals would result
in a better value for $C^*$, but only by a minuscule amount.

Wu gives us that $H(s)$ is decreasing on $[1,10]$. The constraint on the parameter $s$ is that
$2 \leq s \leq 3$, but Wu's discretisation starts at $s =2.2$. 
Wu claims this was chosen as $\Psi_2$ attains its maximum at approximately $s = 2.1$. 
How close this value of $s$ is to the true maximal value is not stated, so this could be another avenue
for improvement. One could explore the value of $\Psi_2$ in a small neighbourhood around $s=2.1$,
to see what the true maximum is.
  
Slightly nudging this value left of $2.2$ to, say,
$2.19$ might mean a better upper bound for $H$, and thus a better value for Chen's constant.

\section{Future work}

There is still a lot of work that could be done to decrease Chen's constant.
The obvious improvement is to find a way to compute $\Psi$\footnote{Whenever $\Psi$ is mentioned in this section,
take it to mean both $\Psi_1$ and $\Psi_2$.} using the approximation
of $\omega(u)$ given above. We could also look at attempting to get the current implementation of 
$\Psi$ running quicker, which would allow sampling more points in the interval $[2,2.5]$,
and finding the coefficients $\phi, s'_i,k_{1,i},k_{2,i},k_{3,i}$ that make $\Psi_2$ maximal, as the current version of
$\Psi_2$ is too slow to optimise over 5 variables simultaneously.
We could look at more clever ways to maximise $\Psi$ using more advanced numerical maximisation methods.
The behaviour of $\Psi$ is not pathological, in the sense that the functions appear to have
some degree of smoothness, and the derivative is not too large. It seems true that $\Psi$ are locally
concave
\footnote{A function $f$ is concave if for any two points $x_1 \leq x_2$, the unique line $g$
that intersects $(x_1,f(x_1))$ and $(x_2,f(x_2))$ satisfies $f(z) \geq g(z)$ for all $x_1 \leq z \leq x_2$.} 
near the local maximum, so once values for the variables have been found that bring $\Psi$ near the
maximum, one could then run a battery of standard convex optimisation techniques on the problem. 

We could also look at alternative forms of approximating $\omega(u)$.
Since we are never concerned with the values of $\omega(u)$ at single points, but rather the integral of $\omega(u)$ over intervals,
we could look to approximations of $\omega_j(u)$ on each interval $[j,j+1]$ that minimise the $L^\infty$ error
$$
||\omega_j(u) - \omega(u)||_{L^\infty [j,j+1]} = \sup_{u \in [j,j+1]} |\omega_j(u) - \omega(u)|,
$$
as Taylor expansions are by no means the best polynomial approximation with respect to $L^\infty$ error.
Chebyshev
\footnote{The Chebyshev polynomials are a set of polynomials defined by the recursive formula
$T_0(x) = 1, T_1(x) = x, T_{n+1}(x) = 2xT_n(x) - T_{n-1}(x)$. They have the property that the error
incurred by using them to approximate a function is more uniformly distributed over the interval the function is
being approximated. Taylor series are exact at a single point, and the error tends to grow quickly as
points move away from the center of the expansion \cite{burden1985numerical}.}
or
Hermite
\footnote{Hermite polynomials are similar to Chebyshev polynomials, but Hermite polynomial approximation
also takes into account the derivative of the function being approximated, and attempts to match both the values of the function,
and it's derivative \cite{burden1985numerical}.}
polynomials are ideal for this purpose, and the speed of computing the approximation is not a concern,
as we only need to compute the polynomial approximation of $\omega(u)$ once. So the new Chebyshev approximation
would be identical to the current approximation, but the polynomials would have different coefficients.
We could also look at Bernstein polynomials.
\footnote{The Bernstein polynomials are a family of polynomials, given by
$B_n(x) = \sum_{k=0}^n \binom{n}{k} x^k(1-x)^{n-k}f(k/n)$. This is the $n^\text{th}$
Bernstein polynomial for the function $f$. It can be shown that 
$||f - B_n||_{L^\infty([0,1])} \to 0$ i.e.\ that $B_n$ converges to $f$ uniformly.\cite{goldberg1964methods}}

 As $\omega(u)$ is continuous, we have by the Weierstrass Approximation Theorem\cite{goldberg1964methods}
that for any $\eps > 0$ and any interval $[1, n]$ there exists a polynomial $p(x)$ such that $||\omega(u) - p(u)||_{L^\infty[1,n]} \leq \eps$.
Bernstein polynomials provide a constructive way to find such a polynomial $p(x)$. This might prove beneficial as integrating 
$\omega(\frac{\phi-t-u-v}{u})$ in $I_1$ (\ref{eq:I1}) 
would prove easier if $\omega(u)$ were represented as a single polynomial instead of a spline, but since $\omega(u)$ is not smooth and in fact has a cusp at $u=1$, the degree of a polynomial required to closely match
$\omega(u)$ to within say $10^6$ everywhere could be massive.

None of these piecewise polynomial implementations of $\omega(u)$ 
would fix the problem of integrating $I_1$ (\ref{eq:I1}), but
if that problem were solved, it might permit Wu's values to be more accurately replicated. 

For more rigorous ways of showing the difficulty of improving $C^*$, we could 
compute some additional points for Wu's table (\ref{tab:wu}), and then interpolate the result.
The interpolation method did remarkably well at demonstrating the negligible improvements to $C^*$ by a finer discretisation.
However, with so few points to interpolate, it is difficult to argue that our interpolant functions for $\Psi$
are actually representative of the behaviour of $\Psi$.
Ideally, we would want to, either by analytic or numeric methods, create an upper bound for $\Psi$ that is in a sense ``too good'',
and then use that to compute $C^*$. If the resulting value for $C^*$ was only a very small change from Wu's value, then
that provides very strong evidence that sizeable improvements on $C^*$ cannot be made by using Wu's method.



  \appendix

  
\chapter{}
\label{appendix1}

All the code used for the results in this thesis was written in \emph{Mathematica}.
Important sections of the code are here, with a brief description of what they relate to.
The full code used in the thesis is available upon request.

\section{Code}
Implementation of $D(N)$

\begin{verbatim}
d[n_] := If[Mod[n/2, 2] == 0, deven[n], dodd[n]]
deven[n_] := Module[{t = 0, p = 2}, 
     While[p <= n/2, If[PrimeQ[n - p], t++]; p = NextPrime[p]]; t]
dodd[n_] := Module[{t = 0, p = 2}, 
     While[p <= n/2 - 1, If[PrimeQ[n - p], t++]; p = NextPrime[p]]; 
  If[PrimeQ[n/2], t++]; t]
\end{verbatim}

Taylor expansion of the Buchstab function $\omega(u)$.
\begin{verbatim}
w2*Exact[u_] := (Log[u - 1] + 1)/u; 
degreeT = 10; 
precisionT = MachinePrecision; 
intervalsT = 5; 
wT[1, u_] := 
  Evaluate[N[Normal[Series[1/u, {u, 1 + 1/2, degreeT}]], precisionT]]; 
wT[2, u_] := 
  Evaluate[N[
    Normal[Series[w2Exact[u], {u, 2 + 1/2, degreeT}]], 
    precisionT]]; 
aw2 = 
  N[Evaluate[
    CoefficientList[wT[2, u] /. u -> x + 2 + 1/2, x]], 
   precisionT]; 
Clear[aw]; 
aw[k_, 2] := aw2[[k + 1]]
aw[0, j_] := 
 aw[0, j] = (1/(j + 1/2))*
   Sum[(aw[k, j - 1]/2^k)*(j + (-1)^k/(2*(k + 1))), 
         {k, 0, degreeT}]
aw[k_, j_] := 
 aw[k, j] = (1/(j + 1/2))*(aw[k - 1, j - 1]/k - aw[k - 1, j])
wT[j_, u_] := 
 ReleaseHold[
  Hold[Sum[aw[k, j]*(u - (j + 1/2))^k, {k, 0, degreeT}]]]
wLim = Exp[-EulerGamma]; 
wA[u_] := Boole[u > intervalsT + 1]*wLim + 
     Sum[Boole[Inequality[k, LessEqual, u, Less, k + 1]]*wT[k, 
     u], {k, 1, intervalsT}]
\end{verbatim}

Numerical solution to $\omega(u)$.
\begin{verbatim}
wN[u_] := Evaluate[w[u] /. 
NDSolve[{D[u*w[u], u] == w[u - 1], w[u /; u <= 2] == 1/u}, w, 
      {u, 1, 30}, WorkingPrecision -> 50]]; 

\end{verbatim}

Defining {\tt unsafeIntegral} to blindly apply FTC without checking for validity.
\begin{verbatim}
unsafeIntegrate[f_] := f;
unsafeIntegrate[f_, var__, 
  rest___] :=
 (unsafeIntegrate[(Integrate[f, var[[1]]] /. 
      var[[1]] -> var[[3]]), rest]
   - unsafeIntegrate[(Integrate[f, var[[1]]] /. var[[1]] -> var[[2]]),
     rest])
\end{verbatim}

Numerically evaluating $I_1$ (\ref{eq:I1})

\begin{verbatim}
i1phiN[p_, sp_, s_] := 
  Quiet[ReleaseHold[NIntegrate[Boole[t <= u <= v]
  	*wN[(p - t - u - v)/u]*(1/(t*u^2*v)), 
      {u, sp, s}, {t, sp, s}, {v, sp, s}, MaxRecursion -> 10, 
      WorkingPrecision -> MachinePrecision, MaxPoints -> 100000, 
      	AccuracyGoal -> 20]]][[1]]
\end{verbatim}

Numerically evaluating $I_{2_i}$ (\ref{eq:I2})
\begin{verbatim}
i2[i_, p_, s_, sp_, k1_, k2_, k3_] := 
  Piecewise[{{NIntegrate[Boole[t <= u <= v]
  	*wN[(p - t - u - v)/u]*(1/(t*u^2*v)), 
       {t, 1/k1, 1/k3}, {u, 1/k1, 1/k3}, {v, 1/k1, 1/k3}, 
       MaxRecursion -> 10, 
       WorkingPrecision -> MachinePrecision, MaxPoints -> 100000, 
       AccuracyGoal -> 20][[1]], 
     i == 9}, 
     {NIntegrate[Boole[t <= u]*wN[(p - t - u - v)/u]*(1/(t*u^2*v)), 
       {t, 1/k1, 1/k2}, {u, 1/k1, 1/k2}, {v, 1/k2, 1/s}, 
       MaxRecursion -> 10, 
       WorkingPrecision -> MachinePrecision, MaxPoints -> 100000, 
       AccuracyGoal -> 20][[1]], 
     i == 10}, 
     {NIntegrate[Boole[u <= v]
     	*wN[(p - t - u - v)/u]*(1/(t*u^2*v)), 
       {t, 1/k1, 1/k2}, {u, 1/k2, 1/k3}, {v, 1/k2, 1/k3},
       MaxRecursion -> 10, 
       WorkingPrecision -> MachinePrecision, MaxPoints -> 100000, 
       AccuracyGoal -> 20][[1]], 
     i == 11}, 
     {NIntegrate[Boole[t <= u]*wN[(p - t - u - v)/u]*(1/(t*u^2*v)), 
       {t, 1/sp, 1/k1}, {u, 1/sp, 1/k1}, {v, 1/k3, 1/s}, 
       MaxRecursion -> 10, 
       WorkingPrecision -> MachinePrecision, MaxPoints -> 100000, 
       AccuracyGoal -> 20][[1]], 
     i == 12}, 
     {NIntegrate[wN[(p - t - u - v)/u]*(1/(t*u^2*v)), 
     {t, 1/sp, 1/k1}, {u, 1/k1, 1/k2}, {v, 1/k2, 1/s}, 
     MaxRecursion -> 10, WorkingPrecision -> 
        MachinePrecision, MaxPoints -> 100000, 
        AccuracyGoal -> 20][[1]], 
        i == 13}, 
    {NIntegrate[Boole[u <= v]*wN[(p - t - u - v)/u]*(1/(t*u^2*v)), 
    {t, 1/sp, 1/k1}, {u, 1/k2, 1/s}, {v, 1/k2, 1/s}, 
    MaxRecursion -> 10, WorkingPrecision -> 
        MachinePrecision, MaxPoints -> 100000, AccuracyGoal -> 20][[1]], 
        i == 14}, 
    {NIntegrate[wN[(p - t - u - v)/u]*(1/(t*u^2*v)), 
    {t, 1/k1, 1/k2}, {u, 1/k2, 1/k3}, {v, 1/k3, 1/s}, 
    MaxRecursion -> 10, WorkingPrecision -> MachinePrecision, 
       MaxPoints -> 100000, AccuracyGoal -> 20][[1]], 
       i == 15}, 
    {NIntegrate[Boole[t <= u <= v <= w]
    *wN[(p - t - u - v - w)/v]*(1/(t*u*v^2*w)), 
       {t, 1/k2, 1/k3}, {u, 1/k2, 1/k3}, {v, 1/k2, 1/k3}, 
       {w, 1/k2, 1/k3}, 
       MaxRecursion -> 10, WorkingPrecision -> MachinePrecision, 
       MaxPoints -> 100000, 
       AccuracyGoal -> 20][[1]], 
       i == 16}, 
    {NIntegrate[Boole[t <= u <= v]
    *wN[(p - t - u - v - w)/v]*(1/(t*u*v^2*w)), 
       {t, 1/k2, 1/k3}, {u, 1/k2, 1/k3}, {v, 1/k2, 1/k3}, 
       {w, 1/k3, 1/s}, 
       MaxRecursion -> 10, WorkingPrecision -> MachinePrecision, 
       MaxPoints -> 100000, 
       AccuracyGoal -> 20][[1]], 
       i == 17}, 
    {NIntegrate[Boole[t <= u && v <= w]
    *wN[(p - t - u - v - w)/v]*(1/(t*u*v^2*w)), 
       {t, 1/k2, 1/k3}, {u, 1/k2, 1/k3}, 
       {v, 1/k3, 1/s}, {w, 1/k3, 1/s}, 
       MaxRecursion -> 10, WorkingPrecision -> MachinePrecision, 
       MaxPoints -> 100000, 
       AccuracyGoal -> 20][[1]], i == 18}, 
    {NIntegrate[Boole[u <= v <= w]
    *wN[(p - t - u - v - w)/v]*(1/(t*u*v^2*w)), 
       {t, 1/k1, 1/k2}, {u, 1/k3, 1/s}, {v, 1/k3, 1/s}, 
       {w, 1/k3, 1/s}, MaxRecursion -> 10, 
       WorkingPrecision -> MachinePrecision, 
       MaxPoints -> 100000, AccuracyGoal -> 20][[1]], 
     i == 19}, 
     {NIntegrate[Boole[u <= v <= w]
     *wN[(p - t - u - v - w - x)/w]*
        (1/(t*u*v*w^2*x)), {t, 1/k2, 1/k3}, {u, 1/k3, 1/s}, 
        {v, 1/k3, 1/s}, {w, 1/k3, 1/s}, 
       {x, 1/k3, 1/s}, MaxRecursion -> 10, 
       WorkingPrecision -> MachinePrecision, 
       MaxPoints -> 100000, 
       AccuracyGoal -> 20][[1]], i == 20}, 
    {NIntegrate[Boole[t <= u <= v <= w <= x <= y]
    	*wN[(p - t - u - v - w - x - y)/x]*
        (1/(t*u*v*w*x^2*y)), {t, 1/k3, 1/s}, 
        {u, 1/k3, 1/s}, {v, 1/k3, 1/s}, 
       {w, 1/k3, 1/s}, {x, 1/k3, 1/s}, {y, 1/k3, 1/s}, 
       MaxRecursion -> 10, 
       WorkingPrecision -> MachinePrecision, MaxPoints -> 100000, 
       AccuracyGoal -> 20][[1]], 
     i == 21}}]

\end{verbatim}

Setting up Wu's data

\begin{verbatim}
k1L[x_] := {3.53, 3.54, 3.57, 3.56}[[x]]; 
k2L[x_] := {2.9, 2.88, 2.87, 2.91}[[x]]; 
k3L[x_] := {2.44, 2.43, 2.4, 2.5}[[x]]; 
sL[i_] := Piecewise[{{1, i == 0}, {2.1 + 0.1*i, i > 0}}]; 
spL[i_] := {4.54, 4.5, 4.46, 4.12, 3.58, 3.47, 3.34, 3.19, 3.}[[i]];
\end{verbatim}

Defining $\Psi_2$

\begin{verbatim}
psi2comp[i_, s_, s2_, k1_, k2_, 
  k3_] := (-2/5)*Integrate[Log[t - 1]/t, {t, 2, s2 - 1}] - 
     (2/5)*Integrate[Log[t - 1]/t, {t, 2, k1 - 1}] - 
     (1/5)*Integrate[Log[t - 1]/t, {t, 2, k2 - 1}] + 
     (1/5)*
   Integrate[Log[s2*t - 1]/(t*(1 - t)), {t, 1 - 1/s, 1 - 1/s2}] + 
     (1/5)*
   Integrate[Log[k1*t - 1]/(t*(1 - t)), {t, 1 - 1/k3, 1 - 1/k1}] + 
     (-2/5)*Sum[i2comp[i, k], {k, 9, 21}]
\end{verbatim}
where {\tt i2comp} gives the output of (Listing \ref{lst:I1})
\begin{verbatim}
i2comp[i_, k_] := Map[Last, \[Phi]Mout[[i]]][[k - 8]]
\end{verbatim}

Part of $\Psi_2$

\begin{verbatim}
i2Sum[i_, s_, sp_, k1_, k2_, k3_] := 
   (-2/5)*Sum[i2[k, \[Phi]Max[k, i], sL[i], spL[i],
   		 k1L[i], k2L[i], k3L[i]], {k, 9, 19}]; 
\end{verbatim}

Implementation of (Listing \ref{lst:I1}), numerical maximisation of $I_1$ with
respect to $\phi$.

\begin{verbatim}
i1testMax[i_, tolmax_] := Module[{tol = 0.1, pMax, pL = pLB[i, k], pH = 5}, 
   While[tol > tolmax, 
     pMax = MaximalBy[monitorParallelTable[Quiet[{p, i2[k, p, sL[i], spL[i], k1L[i], 
              k2L[i], k3L[i]]}], {p, pL, pH, tol}, 1], Last][[1]][[1]]; 
      pL = Max[pMax - tol, 2]; pH = pMax + tol; tol = tol/10]; 
    {pMax, i2[k, pMax, sL[i], spL[i], k1L[i], k2L[i], k3L[i]]}]
\end{verbatim}

Implementation of (Listing \ref{lst:bisect}) to find the smallest $\phi$ such
that the integral $I_{2,i}$ is still defined.

\begin{verbatim}
i1Find[i_, k_, {p_, pL_, pH_, tol_}] := 
  Piecewise[{{2, Quiet[Check[Quiet[i2[k, 2, sL[i], spL[i], k1L[i], k2L[i], 
  k3L[i]]], 
        "Null"]] != "Null"}, {i2Find2[i, k, {p, pL, pH, tol}], True}}]
i2Find[i_, k_, {p_, pL_, pH_, tol_}] := 
  Piecewise[{{(pH + pL)/2, pH - pL < tol && 
      Quiet[Check[Quiet[i2[k, (pL + pH)/2, sL[i], spL[i], k1L[i], k2L[i], k3L[i]]], 
         "Null"]] != "Null"}, {i2Find2[i, k, {p, pL, (pH + pL)/2, tol}], 
     Quiet[Check[Quiet[i2[k, (pL + pH)/2, sL[i], spL[i], k1L[i], k2L[i], k3L[i]]], 
        "Null"]] != "Null"}, {i2Find2[i, k, {p, (pH + pL)/2, pH, tol}], 
     Quiet[Check[i2[k, (pH + pL)/2, sL[i], spL[i], k1L[i], k2L[i], k3L[i]], 
     "Null"]] == 
      "Null"}}]
i2Find[i_, k_, {p_, pL_, pH_, tol_}] := 
  Piecewise[{{2, Quiet[Check[Quiet[i2[k, 2, sL[i], spL[i], k1L[i], k2L[i], 
  k3L[i]]], 
        "Null"]] != "Null"}, {i2Find2[i, k, {p, pL, pH, tol}], True}}]
i2Find2[i_, k_, {p_, pL_, pH_, tol_}] := 
  Piecewise[{{(pH + pL)/2, pH - pL < tol && 
      Quiet[Check[Quiet[i2[k, (pL + pH)/2, sL[i], spL[i], k1L[i], k2L[i], 
      k3L[i]]], 
         "Null"]] != "Null"}, {i2Find2[i, k, {p, pL, (pH + pL)/2, tol}], 
     Quiet[Check[Quiet[i2[k, (pL + pH)/2, sL[i], spL[i], k1L[i], k2L[i], 
     k3L[i]]], 
        "Null"]] != "Null"}, {i2Find2[i, k, {p, (pH + pL)/2, pH, tol}], 
     Quiet[Check[i2[k, (pH + pL)/2, sL[i], spL[i], k1L[i], k2L[i], k3L[i]], 
     "Null"]] == 
      "Null"}}]

\end{verbatim}

Computing $a_{i,j}$ (\ref{eq:aij}).
 \begin{verbatim}
 al[i_, s_, sp_, k1_, k2_, k3_] := Piecewise[{{k1 - 2, i == 1}, 
 {sp - 2, i == 2}, 
      {sp - sp/s - 1, i == 3}, {sp - sp/k2 - 1, i == 4}, 
      {sp - sp/k3 - 1, i == 5}, 
      {sp - (2*sp)/k2, i == 6}, {sp - sp/k1 - sp/k3, i == 7}, 
      {sp - sp/k1 - sp/k2, i == 8}, 
      {k1 - k1/k2 - 1, i == 9}}]; 
 sig0[t_] := sig[3, t + 2, t + 1]/(1 - σ[3, 5, 4]); 
 sig[a_, b_, c_] := 
 	Evaluate[unsafeIntegrate[(1/t)*Log[c/(t - 1)], {t, a, b}]]; 
 Xi1[t_, s_, sp_] := (σ0[t]/(2*t))*Log[16/((s - 1)*(sp - 1))] + 
     (Boole[sp - 2 <= t <= 3]/(2*t))*Log[(t + 1)^2/((s - 1)*(sp - 1))] + 
     (Boole[sp - 2 <= t <= sp - sp/s - 1]/(2*t))*Log[(t + 1)/((s - 1)*
     (sp - 1 - t))]; 
 a1Int[a_, b_, s_, sp_] := NIntegrate[Xi1[t, s, sp], {t, a, b}]; 
 a1[i_, j_] := a1Int[sL[j - 1], sL[j], sL[i], spL[i]]; 

 \end{verbatim}

and $\Xi_2$ (\ref{eq:Xi2})

\begin{verbatim}
Xi2[t_, s_, sp_, k1_, k2_, k3_] := 
   (σ0[t]/(5*t))*Log[1024/((s - 1)*(sp - 1)*(k1 - 1)*(k2 - 1)
   *(k3 - 1))] + 
   (Boole[sp - 2 <= t <= 3]/(5*t))*Log[(t + 1)^5/((s - 1)*
   (sp - 1)*(k1 - 1)*(k2 - 1)*(k3 - 1))] + 
   (Boole[k1 - k1/k2 - 1 <= t <= k1 - 2]/(5*t))*
     Log[(t + 1)/((k2 - 1)*(k1 - 1 - t))] + 
     (Boole[sp - sp/k3 - 1 <= t <= sp - 2]/(5*t))*
     Log[(t + 1)/((k3 - 1)*(sp - 1 - t))] + 
     (Boole[sp - sp/s - 1 <= t <= sp - 2]/(5*t))*
     Log[(t + 1)/((s - 1)*(sp - 1 - t))] + 
     (Boole[k1 - 2 <= t <= sp - 2]/(5*t))*
     Log[(t + 1)^2/((k1 - 1)*(k2 - 1))] + 
     (Boole[sp - sp/k1 - sp/k3 <= t <= sp - sp/k3 - 1]/
      (5*t*(1 - t/sp)))*Log[sp^2/((k1*sp - sp - k1*t)*(k3*sp - sp - k3*t))]
       + (Boole[sp - sp/k3 - 1 <= t <= sp - sp/k1 - sp/k2]/(5*t*(1 - t/sp)))*
     Log[(sp*(sp - 1 - t))/(k1*sp - sp - k1*t)] + 
    (Boole[sp - (2*sp)/k2 <= t <= sp - sp/k1 - sp/k2]/(5*t*(1 - t/sp)))*
     Log[sp/(k2*sp - sp - k2*t)] + (Boole[sp - sp/k1 - sp/k2 <= t <= sp - 2]/
      (5*t*(1 - t/sp)))*Log[sp - 1 - t]; 

a2Int[a_, b_, s_, sp_, k1_, k2_, k3_] 
	:= NIntegrate[Xi2[t, s, sp, k1, k2, k3], {t, a, b}]; 

a2[i_, j_] := a2Int[sL[j - 1], sL[j], sL[i], spL[i], k1L[i], k2L[i], k3L[i]]; 

\end{verbatim}

Solving the matrix equations (\ref{eq:wu_matrix})

\begin{verbatim}
bChen = {0.015826357,
  0.015247971,
  0.013898757,
  0.011776059,
  0.009405211,
  0.006558950,
  0.003536751,
  0.001056651,
  0};
aMat = Table[a[i, j], {i, 1, 9}, {j, 1, 9}];
MatrixForm[LinearSolve[IdentityMatrix[9] - aMat, bVec]]
  
  
\end{verbatim}

Code to interpolate Wu's data and integrate over $d$ many intervals.

\begin{verbatim}
psi1data = Table[{sL[i], ψ1[[i - 4]]}, {i, 5, 9}]; 
psi2data = Table[{sL[i], ψ2[[i]]}, {i, 1, 4}]; 

int1 = Interpolation[psi1data]; 
int2 = Interpolation[psi2data]; 

root = s /. FindRoot[int1[s] - int2[s], {s, 2.55}]; 

interpsi[s_] := 
	Piecewise[{{int2[2.2], 2 <= s <= 2.2}, {int2[s], 2.2 <= s <= root}, 
     {int1[s], root <= s <= 3}}]; 

spdata = Table[{sL[i], spL[i]}, {i, 1, 9}]
k1data = Table[{sL[i], k1L[i]}, {i, 1, 4}]; 
k2data = Table[{sL[i], k2L[i]}, {i, 1, 4}]; 
k3data = Table[{sL[i], k3L[i]}, {i, 1, 4}]; 

intk1 = Interpolation[k1data]; 
intk2 = Interpolation[k2data]; 
intk3 = Interpolation[k3data]; 
intSp = Interpolation[spdata]

a1Int[a_, b_, s_, sp_] := NIntegrate[Xi1[t, s, sp], {t, a, b}]

sint[i_, d_] := sint[i, d] = Piecewise[{{1, i == 0}, 
      {Range[22/10, 3, (3 - 22/10)/(d - 1)][[i]], i > 0}}]; 

a1[i_, j_, d_] 
	:= a1Int[sint[j - 1, d], sint[j, d], sint[i, d], intSp[sint[i, d]]]; 

a2Int[a_, b_, s_, sp_, k1_, k2_, k3_] 
	:= NIntegrate[Xi2[t, s, sp, k1, k2, k3], {t, a, b}]; 

a2[i_, j_, d_] 
	:= a2Int[sint[j - 1, d], sint[j, d], sint[i, d], intSp[sint[i, d]], 
    intk1[sint[i, d]], intk2[sint[i, d]], intk3[sint[i, d]]]; 

a[i_, j_, d_] 
	:= Piecewise[{{a1[i, j, d], sint[i, d] > root}, 
     {a2[i, j, d], sint[i, d] < root}}]; 

aMat[d_] := ParallelTable[a[i, j, d], {i, 1, d}, {j, 1, d}]; 
bsec[i_, d_] := interpsi[sint[i, d]]; 
bVec[d_] := Table[interpsi[sint[i, d]], {i, 1, d}]; 
xVec[d_] := LinearSolve[IdentityMatrix[d] - aMat[d], bVec[d]]; 
chen[d_] := 8*(1 - xVec[d][[1]])

Quiet[Function[x, AbsoluteTiming[{x, chen[x]}]] /@ {200, 300, 400, 500, 600}]
\end{verbatim}

\newpage



  \addcontentsline{toc}{chapter}{Bibliography}

  \bibliography{biblo}
  \bibliographystyle{abbrv}


\end{document}